\documentclass{amsart}
\usepackage{amssymb,hyperref}

 \newcommand{\Z}{{\mathbb Z}}
 \newcommand{\C}{{\mathbb C}}
 
\newcommand{\Q}{{\mathbb Q}}
 
 \newcommand{\T}{{\mathbb T}}

\newcommand{\cal}{\mathcal}

\newcommand{\p}{\partial}

\newcommand{\ve}{\varepsilon}
 \newtheorem{theorem}{Theorem}
 \newtheorem{thm}[theorem]{Theorem}
 \newtheorem{lemma}[theorem]{Lemma}
  \newtheorem{proposition}[theorem]{Proposition}
 \newtheorem{corollary}[theorem]{Corollary}
 \newtheorem{cor}[theorem]{Corollary}
 \newtheorem{remark}[theorem]{Remark}
 \newtheorem{example}[theorem]{Example}

 \newtheorem{question}[theorem]{Question}
\newcommand{\la}{\left\langle}
\newcommand{\ra}{\right\rangle}

\newcommand{\g}{\mathfrak g}

\newcommand{\comment}[1]{}

\title{Character Varieties}
\author{Adam S. Sikora}


\date{}

\begin{document}
\thispagestyle{empty}

\begin{abstract}
We study properties of irreducible and completely reducible representations
of finitely generated groups $\Gamma$ into reductive algebraic groups $G.$
In particular, we study the geometric invariant theory of the $G$ action on the space of $G$-representations of $\Gamma$ by conjugation.

Let $X_G(\Gamma)$ be the $G$-character variety of $\Gamma.$
We prove that for every completely reducible, scheme smooth $\rho:\Gamma\to G$
$$T_{[\rho]}\, X_G(\Gamma)\simeq T_0\, \big(H^1(\Gamma,Ad\, \rho)//S_\Gamma\big),$$
where $H^1(\Gamma,Ad\, \rho)$ is the first cohomology group of $\Gamma$ with coefficients in the Lie algebra $\g$ of $G$ twisted by
$\Gamma\stackrel{\rho}{\longrightarrow} G \stackrel{Ad}{\longrightarrow} GL(\g)$
and $S_\Gamma$ is the centralizer of $\rho(\Gamma)$ in $G.$
The condition of $\rho$ being scheme smooth is very important as there are groups $\Gamma$ such that
$$dim\, T_{[\rho]}\, X_G(\Gamma)< T_0\, H^1(\Gamma,Ad\, \rho),$$
for a Zariski open subset of points in $X_G(\Gamma).$
We prove, however, that all irreducible representations of surface groups are scheme smooth.

Let $M$ be an orientable $3$-manifold with a connected boundary $F$ of genus $g\geq 2.$
Let $X_G^g(F)$ be the subset of the $G$-character variety of $\pi_1(F)$ composed of
conjugacy classes of good representations
$\rho: \Gamma\to G,$ i.e. irreducible representations such that the centralizer of $\rho(\Gamma)$ is the center of $G.$
By a theorem of Goldman, $X_G^g(F)$ is a holomorphic symplectic manifold.
We prove that the set of good $G$-representations of $\pi_1(F)$ which extend to representations of $\pi_1(M)$ is a complex isotropic subspace of $X_G^g(F).$
It is Lagrangian, if these representations correspond to reduced points of the $G$-character variety of $M$. It is an open problem whether it is always the case.
\end{abstract}

\maketitle
\tableofcontents \vspace{.2in}
\pagestyle{myheadings}
\markboth{\hfil{\sc }\hfil}
{\hfil{\sc Character Varieties}\hfil}

\section{Summary of Results}
\label{s_sum}
Let $G$ be a complex reductive algebraic group, for example a classical group of matrices, $GL(n,\C), SL(n,\C), O(n,\C), Sp(n,\C)$ or its quotient.
Let $\Gamma$ be a finitely generated group.
We say that a representation $\rho: \Gamma\to G$ is {\em irreducible} if $\rho(\Gamma)$ is not contained in any proper parabolic subgroup of $G$. Additionally, we say that $\rho$
is {\em completely reducible} if for every parabolic subgroup $P\subset G$ containing $\rho(\Gamma)$,
there is a Levi subgroup $L\subset P$ containing $\rho(\Gamma).$ In particular, $\rho:\Gamma\to GL(n,\C)$ is irreducible if and only if $\C^n$ is a simple $\Gamma$-module (via $\rho$)
and it is completely reducible if $\C^n$ is a semi-simple $\Gamma$-module.

We discuss properties of irreducible and completely reducible representations
in Sections \ref{s_irred}-\ref{s_stab}. For example we prove that $\rho:\Gamma\to G$ is completely reducible if and only if the algebraic closure of $\rho(\Gamma)$ in $G$ is a linearly reductive group. Furthermore, a completely reducible representation $\rho: \Gamma\to G$ is irreducible if and only if the centralizer of $\rho(\Gamma)$ is a finite extension of $C(G).$

The space, $Hom(\Gamma,G),$ of all group homomorphisms from $\Gamma$ to $G$
is an algebraic set on which $G$ acts by conjugation. 
In Section \ref{s_git}, we study properties of this action from the point of view of the Geometric Invariant Theory. In particular we observe that $\rho$ is a poly-stable point under that action if and only if $\rho$ is completely reducible. If $\rho$ is irreducible then it is a stable point. Finally, $\rho$ is properly stable if and only if $\rho$ is irreducible and $C(G)$ is finite.

The categorical quotient, $X_G(\Gamma)=Hom(\Gamma,G)//G,$ is called the {\em $G$-character variety of $\Gamma$}, c.f. Section \ref{s_char_var}. Although it is a coarser quotient than the set theory one, it
has the advantage of having a natural structure of an affine algebraic set. Every element of $X_G(\Gamma)$ is represented by a unique completely reducible representation.
With $Hom(\Gamma,G)$ and $X_G(\Gamma),$ there are naturally associated algebraic schemes ${\cal Hom}(\Gamma,G)$ and
${\cal X}_G(\Gamma)={\cal Hom}(\Gamma,G)//G$ such that $\C[Hom(\Gamma,G)]$ and $\C[X_G(\Gamma)]$ are nil-radical quotients of
${\cal O}({\cal Hom}(\Gamma,G))$ and ${\cal O}({\cal X}_G(\Gamma)),$ c.f. Sections \ref{s_rep_v} and \ref{s_char_var_sch}.

We formulate cohomological descriptions of the tangent spaces to ${\cal X}_G(\Gamma)$ and to $X_G(\Gamma)$ in Section \ref{s_tan_char}. More specifically, we prove the following:
Let $H^*(\Gamma,Ad,\rho)$ denote the group cohomology of $\Gamma$ with coefficients in the Lie algebra $\g$ of $G$ twisted by the homomorphism $$\Gamma\stackrel{\rho}{\longrightarrow} G \stackrel{Ad}{\longrightarrow} End(\g),$$ where $Ad$ is the adjoint representation. Denote the stabilizer of $\rho$ under the $G$ action (by conjugation) on $Hom(\Gamma,G)$ by $S_\rho.$ It is the centralizer of $\rho(\Gamma)$ in $G.$
There is a natural action of $S_\rho$ on $H^1(\Gamma,Ad\, \rho),$ c.f. Sec. \ref{s_tan_char}.

We call an irreducible $\rho:\Gamma\to G$ {\em good} if the stabilizer of its image, $S_\rho$ is the center of $G$. We say that $\rho$ is {\em scheme smooth} if ${\cal Hom}(\Gamma,G)$ is non-singular at $\rho.$

\begin{theorem}\label{i_tang}(c.f. Theorem \ref{tangent_X})\\
(1) For every good $\rho:\Gamma\to G$ there exists a natural linear isomorphism
 $$\phi: H^1(\Gamma,Ad\, \rho)\hookrightarrow T_{[\rho]}\, {\cal X}_G(\Gamma).$$
(2) If $\rho:\Gamma\to G$ is completely reducible and scheme smooth then
$$T_0\, \left( H^1(\Gamma,Ad\, \rho)//S_\rho \right)\simeq T_{[\rho]}\, {\cal X}_G(\Gamma)=T_{[\rho]}\, X_G(\Gamma).$$
\end{theorem}

The quotient on the left side in the statement of Theorem \ref{i_tang}(2) may be non-trivial even if $\rho$ is irreducible, c.f. discussion in Section \ref{s_stab} and Example \ref{example-Klein} in particular. We discuss the question of the existence of a natural isomorphism in part (2) of the above theorem in Sec. \ref{s_tan_char}.

Versions of the above theorem for $G=PSL(2,\C)$ appear in \cite[Prop 3.5]{Po} and \cite[Prop. 5.2]{HP2}.
Although statements similar to Theorem \ref{i_tang}(1) appear in literature for general $G$, there is a widespread confusion about the necessary assumptions and all proofs known to us are incomplete.

The following result illuminates the importance of the requirement of $\rho$ being scheme smooth in Theorem \ref{i_tang}(2). Denote the set of equivalence classes of good representations in ${\cal X}_G(\Gamma)$ by ${\cal X}^g_G(\Gamma).$

\begin{theorem}(c.f. Theorem \ref{H^1_iff_reduced})\\
${\cal X}_G^g(\Gamma)$ is reduced iff
$T_{[\rho]}\, X_G(\Gamma)= H^1(\Gamma,Ad\, \rho)$ for a non-empty Zariski open set of $[\rho]$'s in $X_G^g(\Gamma)$.
\end{theorem}

It is easy to see that all representations of free groups are scheme smooth. We prove
an analogous statement for surface groups, i.e. fundamental groups of closed, orientable surfaces of genus $\geq 2$.

\begin{theorem}\label{i_reduced}
For every reductive $G$, all irreducible $G$-representations of surface groups are scheme smooth.
\end{theorem}

This theorem has some important consequences to the theory of skein modules, c.f. Corollary \ref{skein}:

\begin{theorem}
For every orientable surface $F$, there is an isomorphism between
the Kauffman bracket skein algebra ${\cal S}_{2,\infty}(F,\C,-1)$ and
$\C[X_{SL(2,\C)}(\pi_1(F))]$ sending a link $K_1\cup ...\cup K_n$ in $F\times [0,1]$
to $(-1)^n\tau_{K_1}\cdot ...\cdot \tau_{K_n},$ where $\tau_K\in \C[X_{SL(2,\C)}(\pi_1(F))]$ sends the equivalence class of $\rho: \pi_1(F)\to SL(2,\C)$
to $tr\rho(\gamma)$ and $\gamma\in \pi_1(F)$ is any element representing $K$ (with some orientation).
\end{theorem}

Although this result was announced in \cite[Thm. 7.3]{PS2}, its proof relied on \cite[Thm 7.4]{PS2} whose proof was not provided. An alternative proof of the above statement was provided independently by L. Charles and J. March\'e in \cite{CM}.


Denote the set of irreducible $G$-representations of $\Gamma$ by $Hom^i(\Gamma,G).$
It is a Zariski open subset of $Hom(\Gamma,G)$, c.f. Proposition \ref{irred_open}.
Since each equivalence class in $Hom(\Gamma,G)//G$ contains a unique closed orbit and the orbit of
every irreducible representation is closed (Proposition \ref{tfae-closed}), the categorical quotient of $Hom(\Gamma,G)$ restricted to $Hom^{i}(\Gamma,G)$ is the set theoretic quotient. Denote $Hom^{i}(\Gamma,G)/G$
by $X_G^{i}(\Gamma).$

\begin{proposition}\label{i_mfld}
Let $\Gamma$ be a free group or a surface group. Then\\
(1) $X_G^i(\Gamma)$ is an orbifold.\\
(2) If $G=GL(n,\C)$ or $SL(n,\C)$ then $X_G^i(\Gamma)$ is a manifold. (See also \cite{FL2}.)
\end{proposition}

We do not know if Proposition \ref{i_mfld}(2) holds for any reductive groups other than $GL(n,\C)$ and $SL(n,\C)$, c.f. Question \ref{question_CI} and Proposition \ref{orbifold}.


We say that $\rho$ is {\em good} if $\rho$ is irreducible and $S_\rho$ is the center of $G$. Denote the set of all such representations by $Hom^{g}(\Gamma,G).$
It is a Zariski open subset of $Hom(\Gamma,G)$, c.f. Proposition \ref{g_open}.
Let $$X_G^{g}(\Gamma)=Hom^{g}(\Gamma,G)//G=Hom^{g}(\Gamma,G)/G.$$
By the above discussion $X_G^{g}(\Gamma)$ is an open subset of $X_G^i(\Gamma)$. Furthermore, it is a smooth manifold for free groups and surface groups $\Gamma.$

For a topological space $Y,$ we abbreviate $X_G(\pi_1(Y))$ by $X_G(Y).$
Let $F$ be a closed orientable surface of genus $\geq 2.$
Goldman proved that every non-degenerate symmetric bilinear $Ad$-invariant form $B$ on
$\g$ gives rise to a holomorphic symplectic $2$-form $\omega_B$ on $X_G^{g}(F)$, \cite{Go-symp}, c.f. Section \ref{s_char_var_surf}.

According to folk knowledge, if $M$ is a compact orientable $3$-manifold with a connected boundary $F$ then the image of the map $r^*: X_G(M)\to X_G(F)$ induced by the embedding $r:F=\p M\hookrightarrow M$ is ``roughly speaking" a Lagrangian subspace of the $G$-character variety of $F$.
In Section \ref{s_lagr}, we formulate this claim precisely and we prove
it in detail. Let $Y_G(M)$ be the non-singular part of $$X^g_G(F)\cap r^*(X_G(M)).$$

\begin{theorem}\label{i-Lagrangian}
(1) $Y_G(M)$ is an isotropic submanifold of $X_G^g(F)$ with respect to $\omega_B.$\\
In particular, every connected component of $Y_G(M)$ of dimension $\frac{1}{2} dim\, X_G(F)$ is Lagrangian.\\
(2) If a connected component $C$ of $Y_G(M)$ contains the conjugacy class of a
representation which is a reduced point of ${\cal Hom}(\pi_1(M),G)$, then $C$ is a Lagrangian submanifold of $Y_G(M).$
\end{theorem}


Although we do not know any example of a $3$-manifold $M$ such that $Y_G(M)$ is not Lagrangian, in light of the above theorem and the results of M. Kapovich mentioned in Section \ref{s_char_var_sch} we believe that such examples exist.

Here is a different version of (2) above:

\begin{theorem}\label{i-Lagrangian2}
If $X_G^s(M)$ denotes the set of equivalence classes of scheme smooth representations in $X_G(M)$
then $X^g_G(F)\cap r^*(X_G^s(M))$ is an immersed Lagrangian submanifold of $X^g_G(F).$
\end{theorem}

Note however that $X_G^s(M)$ may be empty even if $\pi_1(M)$ has good $G$-representations.

Theorems \ref{i-Lagrangian} and \ref{i-Lagrangian2} are relevant to Chern-Simons theory, c.f. Section \ref{s_lagr}.

In the paper we assume familiarity with basic algebraic geometry and the theory of algebraic groups.
The standard references for these topics are \cite{Ha,Shf,Bo,Hu}.

We would like to thank S. Baseilhac, H. Boden, F. Bonahon, W. Goldman, C. Frohman, M. Kapovich, S. Lawton, J. Porti, and the referee for helpful comments.

\section{Reductive Groups}
\label{s_reductive}

Every algebraic group $G$ contains a unique maximal normal connected solvable subgroup, called its {\em radical} and denoted by $Rad\, G.$ A connected group $G$ is {\em semi-simple} if and only if $Rad\, G$ is trivial.
A connected group $G$ is {\em reductive} if and only if $Rad\, G$ is an algebraic torus, $(\C^*)^n.$
In particular, $\C^*$ and all classical matrix groups, $SL(n,\C),$ $O(n,\C)$, $Sp(n,\C)$, are reductive. Furthermore, Cartesian products, quotients and finite connected covers of reductive groups are reductive.
In fact, all reductive groups can be obtained in this way from simple algebraic groups.

Denote the center of $G$ by $C(G)$ and the connected component of the identity in $C(G)$ by $C^0(G).$
For every reductive $G$, $C^0(G)=(\C^*)^n.$
A reductive group $G$ is semi-simple if and only if $C^0(G)$ is trivial.

Let $[G,G]$ be the commutator of $G$. If $G$ is reductive then $[G,G]$ is semi-simple. Furthermore, by \cite[Proposition IV.14.2]{Bo}, the epimorphism
\begin{equation}\label{nu}
\nu: C^0(G)\times [G,G] \to G,\quad \nu(g,h)=g\cdot h
\end{equation}
has a finite kernel, isomorphic to
$$K= C^0(G)\cap [G,G].$$
Therefore, there is a finite quotient
\begin{equation}\label{pi}
\pi: G\to C^0(G)/K\times [G,G]/K.
\end{equation}

An algebraic group $G$ is {\em linearly reductive} if its all $GL(n,\C)$-representations are completely reducible. $G$ is linearly reductive if and only if the connected component, $G^0,$ of its identity is reductive. (This property does not hold for groups over fields of non-zero characteristic.) Therefore, linearly reductive groups are ``virtually reductive".

A maximal connected solvable subgroup of $G$ is called a {\em Borel subgroup}.
A closed subgroup $P\subset G$ is {\em parabolic} if one of the following equivalent
conditions holds: (a) $G/P$ is a complete variety, (b) $G/P$ is a projective variety,
(c) $P$ contains a Borel subgroup of $G,$ c.f. \cite{Bo}.

A Levi subgroup of an algebraic group $H$ is a connected subgroup $L\subset H$
such that $H$ is a semi-direct product of $L$ and the unipotent radical of $H.$
Since $L$ is isomorphic to the quotient of $H$ by its unipotent radical, it is always reductive.
By a result of Mostow, every algebraic group contains a Levi subgroup, c.f. \cite[IV.11.22]{Bo}

\section{Irreducible and Completely Reducible Subgroups}
\label{s_irred}

We say that a subgroup $H$ (closed or not) of $G$ is {\em irreducible} if it is not contained in any proper parabolic subgroup of $G.$ We also say that $H$ is {\em completely reducible} if for every parabolic subgroup $P\subset G$ containing $H$, there is a Levi subgroup of $P$ containing $H$ as well, \cite{Se,BMR}. In particular, every irreducible subgroup is completely reducible.

The following is an important characterization of completely reducible subgroups:

\begin{proposition}\label{cr=r}
For every reductive $G$, a subgroup $H\subset G$ is completely reducible if and only if the algebraic closure of $H$ in $G$ is a linearly reductive group.
\end{proposition}

\begin{proof}
$\Rightarrow$
(1) Assume first that $H$ is irreducible, i.e. not contained in any proper parabolic subgroup of $G.$
Let $\overline H$ be the Zariski closure of $H$ in $G$ and $Rad_u(\overline H)$ be the unipotent radical of $\overline H.$ Let $P={\mathcal P}(Rad_u(\overline H))$ be the parabolic subgroup defined in \cite[30.3]{Hu}. Then $\overline H\subset N_G(Rad_u(\overline H))$ and $N_G(Rad_u(\overline H))\subset P$ by \cite[30.3 Corollary A]{Hu}.
If $\overline H$ is not linearly reductive then $Rad_u(\overline H)$ is non-trivial and
$Rad_u(\overline H)\subset Rad_u(P)$ by \cite[30.3 Corollary A]{Hu}. Therefore, $P$ is a proper subgroup of $G.$

(2) Now we carry the proof in full generality by induction with respect to $dim\, G:$
If $dim\, G=1$ then $G=\C^*$ and the statement holds.
Assume now that it holds for all reductive algebraic groups $G$ of dimension less than $n$.
Let $dim\, G=n.$ If $H$ lies in a proper parabolic subgroup of $G$
then it also lies in a Levi subgroup of $P$ and the statement follows from inductive hypothesis.
If $H$ does not lie in a proper parabolic subgroup of $G$ then $H$ is irreducible in $G$ and the statement follows from part (1).\\

$\Leftarrow$ Suppose $\overline H$ is linearly reductive and $H\subset P.$ Since $P$ is closed,
$\overline H\subset P.$ Now the statement follows from the fact that
every closed linearly reductive subgroup of $P$ lies in a Levi subgroup of $P$.
Since we do not know a good reference to this fact, we enclose its proof here: There is an exact sequence
$$\{e\}\to Rad_U\, P\to P\stackrel{\tau}{\longrightarrow} L \to \{e\},$$
where $Rad_U\, P$ is the unipotent radical of $P.$
Since ${\overline H}^0$ is reductive, it has no connected unipotent subgroups and, therefore,
$\tau$ is an embedding of ${\overline H}^0$ into $L.$ Therefore, the kernel $K$ of $\tau$ restricted ${\overline H}$ is finite. By \cite[Corollary 4.8]{Bo}, $Rad_U\, P$ is a subgroup of upper triangular matrices and, therefore, it has no elements of finite order. Hence, $K$ is trivial.
\end{proof}

A representation $\phi: \Gamma \to G$ is {\em irreducible} or {\em completely reducible} if $\phi(\Gamma)\subset G$ is.
In particular, a representation $\rho:\Gamma\to GL(n,\C)$ is irreducible if and only if
$\C^n$ does not have any $\Gamma$-invariant subspaces other than $\{0\}$ and $\C^n.$
Additionally, $\rho: \Gamma\to GL(n,\C)$ is completely reducible if and only if $\C^n$ decomposes into a sum of irreducible $\Gamma$-modules.

Since a quotient of a reductive group is reductive, Proposition \ref{cr=r} implies:

\begin{corollary}\label{cr_cr}
For every homomorphism $\phi: G_1\to G_2$ of reductive groups the image of completely reducible subgroup of $G_1$ is completely reducible in $G_2.$
\end{corollary}

Similarly, we have:

\begin{lemma}\label{ir_epi}
For every epimorphism $\phi: G_1\to G_2$ of reductive groups, the image of an irreducible subgroup of $G_1$ is irreducible in $G_2.$
\end{lemma}

\begin{proof} Suppose that $\phi(H)$ lies inside a proper parabolic subgroup $P\subset G_2.$
Since $\phi$ induces an isomorphism $G_1/\phi^{-1}(P)\to G_2/P$ and
$G_2/P$ is complete, $G_1/\phi^{-1}(P)$ is complete as well, implying that $\phi^{-1}(P)$
is a proper parabolic subgroup of $G_1$ containing $H.$
\end{proof}

The following example shows that irreducibility of $H \subset G_1$ does not imply irreducibility of $\phi(H)\subset G_2$ if $\phi: G_1\to G_2$ is not an epimorphism, even if $\phi$ is irreducible itself. 

\begin{proposition}\label{not_absolutely}
Let $H=\{A: A\cdot A^T=\pm I\}\subset SL(2,\C).$\\
(1) $H$ is isomorphic to $O(2,\C).$\\
(2) $H\subset SL(2,\C)$ is irreducible.\\
(3) The image of $H$ under the adjoint representation $Ad: SL(2,\C)\to SL(3,\C)$ is completely reducible but not irreducible in $SL(3,\C).$
\end{proposition}

\begin{proof}
(1) 
$H$ is a non-abelian split $\Z/2$ extension of $SO(2)\simeq \C^*.$
However, since $H^2(\Z/2, \C^*)=0$ (with twisted coefficients), $O(2,\C)$ is the unique non-abelian split extension of $SO(2,\C)$ by $\Z/2$.

(2)  Since $H$ is reductive, it is completely reducible in $SL(2,\C)$ by Proposition \ref{cr=r}. If it was reducible, it would be a subgroup of diagonal matrices, $\C^*.$
Since $H$ is nonabelian, it is irreducible in $SL(2,\C).$

(3) Complete reducibility follows from Corollary \ref{cr_cr}. We claim that the group $Ad(H)$ lies in the parabolic subgroup of
$SL(3,\C)$ composed of transformations of $sl(2,\C)\simeq \C^3$ which preserve $Span(M)\subset sl(2,\C),$
where $M=\left(\begin{matrix} 0 & -1\\ 1 & 0\\ \end{matrix}\right).$
Indeed, since $SL(2,\C)=Sp(2,\C),$ $AMA^T=M$ for every $A\in SL(2,\C).$
If $A\in H$ then $A^T=\pm A^{-1}$ and the claim follows.
\end{proof}

We say that $H\subset G$ is {\em $Ad$-irreducible} if $Ad(H)\subset GL(\g)$ is irreducible. The $H\subset SL(2,\C)$ above is irreducible but not $Ad$-irreducible.
By Proposition \ref{cr=r}, every irreducible subgroup $H\subset G$ is {\em completely $Ad$-reducible}, i.e. $Ad(H)\subset GL(\g)$ is completely reducible.
We are going to show that $Ad$-irreducibility implies irreducibility.

\begin{lemma}\label{irred_irred}
Let $\phi: G\to GL(n,\C)$ be an irreducible representation of a reductive group $G$. If $H$ a subgroup of $G$ such that $\phi(H)$ is irreducible then either\\
(a) $H$ is irreducible, or\\
(b) $Ker\, \phi$ contains the unipotent radical (i.e. the maximal connected unipotent subgroup) of a Borel subgroup of $G.$
\end{lemma}

\begin{proof}
Suppose that $H\subset P\subsetneq G$. Then $\phi$ restricted to $P$ is irreducible as well. Let $U$ be the unipotent radical of $P.$
Denote the space of vectors in $V=\C^n$ invariant under the $U$-action by $V^{U}.$
Since $P$ is a semi-direct product of $U$ and its Levi subgroup $L,$ \cite{Bo}, $P=L \rtimes U$,
$l^{-1}ul\in U$, for every $u\in U$ and $l\in L$,
and
$$u\cdot l\cdot v=l\cdot l^{-1}ul\cdot v=l\cdot v\quad \text{for every $v\in V^{U}$}.$$
Therefore $l\cdot v\in V^{U}$ and, consequently, $V^{U}$ is preserved by $P.$
Since $\phi$ restricted to $P$ is irreducible, by Shur's Lemma $V^{U}$ is either $0$ or $V.$
However, $U$ is a connected solvable group and, therefore, $V^{U}\ne 0$, by Lie-Kolchin theorem,
\cite[Cor 10.5]{Bo}. Hence $V^{U}=V$ and, consequently, $U\subset Ker\, \phi.$ If $B$ is a Borel subgroup of $G$ contained in $P$ then the unipotent radical of $B$ is contained in $U.$
\end{proof}

Since the kernel of the adjoint representation is the center of $G$, \cite[I.3.15]{Bo}, and its
unipotent radical is trivial, Lemma \ref{irred_irred} implies:

\begin{corollary}\label{Ad-cor}
Every $Ad$-irreducible subgroup of a reductive group is irreducible.
\end{corollary}

%
\section{Stabilizers of irreducible representations}
\label{s_stab}
%

\begin{proposition}\label{center-ad}
The centralizer of an $Ad$-irreducible subgroup of a reductive group $G$ is the center of $G.$
\end{proposition}

\begin{proof}
Let $H\subset G$ be $Ad$-irreducible.
By Shur's Lemma the centralizer of $Ad(H)$ is the group of scalar matrices in $GL(\g).$
Hence, $$Ad(C_G(H))\subset C_{GL(\g)}(Ad(H))\subset \{c\cdot I: c\in \C^*\}.$$
On the other hand, since the center of $Ad(G)=G/C(G)$ is trivial, c.f. \cite[Thm 23.16]{FH},
$Ad(G)\cap \{c\cdot I: c\in \C^*\}=\{I\}.$ Hence, $Ad(C_G(H))$ is trivial, implying that
$C_G(H)\subset C(G).$
\end{proof}

\begin{proposition}\label{center}
The centralizer of an irreducible subgroup of a reductive group $G$ is a finite extension of the center of $G.$
\end{proposition}

\begin{proof}
Suppose that the centralizer, $C_G(H),$ of $H$ is an infinite extension of the center. Let $T$ be a maximal
torus in $C_G(H).$ Then $rank\, T>rank\, C(G)$ and $H\subset C_G(T).$ Recall that $T$ is either regular,
semi-regular, or singular, \cite[\S 13.1]{Bo}. If $T$ is regular then $C_G(T)$ is a maximal torus.
If $T$ is semi-regular, then $C_G(T)$ is contained in a Borel subgroup, c.f. proof of \cite[IV.13.1 Proposition]{Bo}. In either case $H\subset C_G(T)$ is reducible. Therefore, $T$ is singular.
In that case $T$ is the connected component of identity of $\bigcap_{\alpha\in I} Ker\, \alpha,$ where the intersection is over a certain proper, non-empty subset $I$ of positive roots. By \cite[IV.14.17]{Bo}, $T$ lies inside of a proper parabolic subgroup of $G$ (denoted by Borel by $P_I$).
\end{proof}

The following lemma will be useful later.

\begin{lemma}\label{centralizer}
For every Levi subgroup $L$ of every proper parabolic subgroup of a reductive group $G$,
$dim\, C(L)>dim\, C(G).$
\end{lemma}

\begin{proof}
Suppose that $L$ is a Levi subgroup of a parabolic group $P$ in $G.$
As before, let $C^0(G)$ be the connected component of the identity in the center of $G.$
Then $C^0(G)\subset L$ and $L'=L/C^0(G)$ is a Levi subgroup of the parabolic subgroup $P'=P/C^0(G)$
of $G'=G/C^0(G).$
Therefore, without loss of generality we can assume that $dim\, C(G)=0$.

By \cite[Prop. 11.23]{Bo}, it is enough to prove that the radical
of $P$, ${\mathcal R}P,$ has a positive dimension.
Fix a root system for $G$. We are going to use the notation of \cite{Bo}.
By classification of parabolic subgroups in \cite[\S 14.17]{Bo}, $P=P_I$ for some subset $I$ of positive roots $\Delta$ of $G.$ Let $T_I$ be the identity component of $\cap_{\alpha\in I} Ker\, \alpha$.
By \cite[Prop. 14.18]{Bo}, $T_I\subset {\mathcal R}P.$ Since $T_I$ is an algebraic torus of dimension
\mbox{$rank\, G-|I|,$} $dim\, {\mathcal R}P>0$ unless $I=\Delta.$ In this case, $P=G.$
\end{proof}

Proposition \ref{center} and Lemma \ref{centralizer} imply:

\begin{corollary}\label{cr_irr} A completely reducible subgroup $H\subset G$ is irreducible if and only if
$dim\, C_G(H)=dim\, C(G).$
\end{corollary}

We will say that a reductive group $G$ has property CI if the centralizer of every irreducible
subgroup of $G$ coincides with the center of $G$.

\begin{example}\label{eg-gln-sln}
$G=GL(n,\C)$ and $SL(n,\C)$ are CI. Indeed, $H\subset G$ is irreducible if and only if elements of $H$
linearly span $M(n,\C).$ Consequently, the centralizer of every irreducible subgroup $H\subset G$ is the center of $G$.
\end{example}

\begin{question}\label{question_CI}
Are $GL(n,\C)$ and $SL(n,\C)$ the only CI-groups?
\end{question}

\begin{example}\label{example-Klein}
$PSL(2,\C)$ is not CI. To see that consider the subgroup $H\subset PSL(2,\C)$, generated by $\pm\left(\begin{matrix} i & 0 \\ 0 & -i \end{matrix}\right)$ and $\pm \left(\begin{matrix} 0 & 1 \\ -1 & 0 \end{matrix}\right).$
Since $|H|=4$, $H$ is the Klein group. One can easily see that $H$ is its own centralizer in $PSL(2,\C)$ (while the center of $PSL(2,\C)$ is trivial).
Being finite, $H$ is linearly reductive and completely reducible by Proposition \ref{cr=r}.
By Corollary \ref{cr_irr}, $H$ is irreducible in $PSL(2,\C).$
\end{example}

\begin{example}
$SO(n,\C)$ is not CI:
Let $DM_n$ be the group of diagonal matrices in $$SO(n,\C)=\{A: A\cdot A^T=I\}\subset SL(n,\C).$$
Then $DM_n\simeq (\Z/2)^{n-1}$ and it is easy to see that $DM_n$ is its own centralizer in $SO(n,\C).$
Being finite, $DM_n$ is linearly reductive and completely reducible by Proposition \ref{cr=r}.
By Corollary \ref{cr_irr}, $DM_n$ is irreducible in $SO(n,\C).$
\end{example}

\begin{proposition}
$Sp(2n,\C)$ is not CI.
\end{proposition}

\begin{proof}(based on the idea of S. Lawton, c.f. \cite{FL2})
Denote by $D(\alpha_1,...\alpha_n)$ the diagonal matrix with entries
$\alpha_1,...,\alpha_n,$ and by $AD(\alpha_1,...\alpha_n)$ the anti-diagonal matrix
$$\left(
\begin{array}{cccc}
 0 & 0 & 0 & \alpha_1 \\
 0 & 0 & ... & 0 \\
 0 & ... & 0 & 0 \\
 \alpha_n & 0 & 0 & 0
\end{array}
\right).$$
The matrices $D(\alpha_1,...,\alpha_n,\alpha_n^{-1},...,\alpha_1^{-1})$ and $AD(\beta_1,....,\beta_n,-\beta_n^{-1},...,-\beta_1^{-1})$, for
$\alpha_1,...,\alpha_n,\beta_1,...,\beta_n\in \C^*$ form a subgroup of
$Sp(2n,\C)=\{A: A J A^T=J\},$ where $J=AD(1,...,1,-1,...,-1).$
Denote that subgroup by $H_n.$ An elementary computation shows that the center of
$H_n$ is composed of matrices $D(\alpha_1,...,\alpha_n,\alpha_n^{-1},...,\alpha_1^{-1}),$
where $\alpha_1,...,\alpha_n\in \{\pm 1\}.$ Since $H_n$ is a finite extension of $(\C^*)^n$,
it is linearly reductive and, hence, by Proposition \ref{cr=r}, it is completely reducible in $Sp(2n,\C).$
Since $C(\Gamma_n)$ is a finite extension of
$C(Sp(2n,\C))=\{\pm 1\},$ $\Gamma_n$ is irreducible by Corollary \ref{cr_irr}.
\end{proof}

By the following result, $PSO(n,\C), PSp(2n,\C)$ are not CI either.

\begin{proposition}
A quotient of a non-CI group by a finite subgroup is non-CI.
\end{proposition}

\begin{proof}
Let $\Gamma\subset G$ be irreducible and such that $C_G(\Gamma)$ is a proper extension of $C(G).$ If $\pi: G\to G'$ is a quotient with finite kernel then $Ker\, \pi\subset C(G)$ and, consequently, the centralizer of $\pi(\Gamma)$ in $G'$ is a proper extension of $C(G').$ Now the statement follows from Proposition \ref{ir_epi}.
\end{proof}

\section{Representation Varieties}
\label{s_rep_v}

If $\Gamma$ is a finitely generated group and $G$ an affine complex algebraic group, then
the space of all $G$-representations of $\Gamma,$ $Hom(\Gamma,G),$ is an algebraic set.

\begin{remark}\label{free_prod}
$Hom(\Gamma_1*\Gamma_2,G)=Hom(\Gamma_1,G)\times Hom(\Gamma_2,G).$
Hence, for the free group on $n$ generators, $Hom(F_n,G)=G^n.$
\end{remark}

\begin{example}\label{Z^2}
Each point of $Hom(\Z^2,SL(2,\C))$ is represented by $\rho: \Z^2\to SL(2,\C)$
defined by
$$\rho(1,0)=\left(\begin{matrix} x_1 & x_2\\ x_3 & x_4\\ \end{matrix}\right),\quad \rho(0,1)=\left(\begin{matrix} x_5 & x_6\\ x_7 & x_8\\ \end{matrix}\right),$$
satisfying relations
$$x_1x_4-x_2x_3-1\ =\ x_5x_8-x_6x_7-1\ =\ x_2x_7-x_3x_6\ =0,$$
$$-x_2x_5+x_1x_6-x_4x_6+x_2x_8\ =\ x_3x_5-x_1x_7+x_4x_7-x_3x_8\ =\ 0.$$
The algebraic set $Hom(\Z^2,SL(2,\C))\subset \C^8$ is irreducible by
\cite[Thm C]{Ri1}.
\end{example}

For a more through study of representation varieties it is useful
to associate with each $\Gamma$ and $G$ as above an affine algebraic scheme, also called the representation variety, whose set of closed points coincides with $Hom(\Gamma,G).$
That scheme, containing sometimes more subtle information about $G$-representations of $\Gamma$ than
$Hom(\Gamma,G),$ is constructed below.

If $G$ is an affine complex algebraic group, then $\C[G]$ is a Hopf algebra with the coproduct
$$\Delta: \C[G]\to \C[G]\otimes \C[G]=\C[G\times G]$$
being the dual to the group product $G\times G\to G$
and the antipode $$S: \C[G]\to \C[G]$$
being the dual to the inverse map $g\to g^{-1}.$
Consequently, for any commutative $\C$-algebra $A$ with product $m: A\times A\to A,$
the space of algebra homomorphisms, $Hom(\C[G],A),$ is a group with
the multiplication
$$Hom(\C[G],A)\times Hom(\C[G],A)\ni (f,g) \to m (f\otimes g)\Delta \in Hom(\C[G],A)$$
and the inverse
$$Hom(\C[G],A)\ni f \to fS \in Hom(\C[G],A).$$
We will denote $Hom(\C[G],A)$ with that group structure by $G(A).$
The functor $G(\cdot)$ is an affine group scheme, \cite{Wa}.
For example, $G(A)=SL(n,A)$ for $G=SL(n,\C).$

We say that a commutative $\C$-algebra $R(\Gamma,G)$ is a {\em universal representation algebra of $\Gamma$ into $G$} and $\rho_U: \Gamma\to G(R(\Gamma,G))$
is a {\em universal representation} if for every commutative
$\C$-algebra $A$ and every representation $\rho: \Gamma\to G(A),$ there is a $\C$-algebra homomorphism $f:R(\Gamma,G)\to A$ inducing a representation $G(f): G(R(\Gamma,G))\to G(A)$ such that
$\rho=G(f)\rho_U,$ \cite{BH,Si1}.

\begin{lemma}\label{universal_prop}
For every $\Gamma$ and every $G$ as above,\\
(1) $R(\Gamma,G)$ and $\rho_U$ exist.\\
(2) $R(\Gamma,G)$ is well defined up to an isomorphism of $\C$-algebras.\\
(3) $\rho_U:\Gamma\to G(R(\Gamma,G))$ is unique up to a composition with $G(f)$
where $f$ is a $\C$-algebra automorphism of $R(\Gamma,G).$\\
\end{lemma}

\begin{proof}
(1) Since each affine algebraic group $G$ is a closed subgroup of $GL(n,\C)$, the coordinate ring
of $G$ is a quotient of
$$\C[GL(n,\C)]=\C[d, x_{ij},\ 1\leq i,j\leq n]/(d\cdot det(x_{ij})-1).$$
Let $$\C[G]=\C[d,x_{ij},\ 1\leq i,j\leq n]/I_G,$$
for an appropriate ideal $I_G$.
For the free group, $F_N=\la \gamma_1,...,\gamma_N\ra$,
$$R(F_N,G)=\C[d_1,x_{1ij},\ 1\leq i,j\leq n]/I_G\otimes ...\otimes \C[d_N,x_{Nij},\
1\leq i,j\leq n]/I_G$$
and
\begin{equation}\label{rho_U}
\rho_U(\gamma_t)=(x_{tij})\in G(R(\Gamma,G)),\quad \text{for}\ t=1,...,N
\end{equation}
satisfy the required universal properties.

If $$\Gamma=\la \gamma_1,...,\gamma_N\ra/H,$$
then we define $R(\Gamma,G)$ as the quotient of $R(\la \gamma_1,...,\gamma_N \ra,G)$
by an ideal $I$ generated by all relations necessary for
(\ref{rho_U}) to be a well defined group homomorphism.
Therefore, each normal generator of $H\triangleleft \la \gamma_1,...,\gamma_N \ra$
introduces $n^2$ relations to $I$ (although some of them may be redundant).
It is easy to see that, in this way, (\ref{rho_U}) descends to a universal representation
$\rho_U: \Gamma \to G(R(\Gamma,G)).$

(2) and (3) follow immediately from the definition.
\end{proof}

Every $\rho\in Hom(\Gamma,G)$ defines a $\C$-algebra homomorphism
$h_\rho: R(\Gamma,G)\to \C$
(unique up to an automorphism of $R(\Gamma,G)$) such that
$$\rho=G(h_\rho)\rho_U.$$
$Ker\, h_\rho$ is a maximal ideal in $R(\Gamma,G)$
and, hence, a closed point in the affine scheme ${\cal Hom}(\Gamma,G)=Spec\, R(\Gamma,G).$
Conversely, every closed point in ${\cal Hom}(\Gamma,G)$ defines a representation $\rho:\Gamma\to G.$ Therefore, $Hom(\Gamma,G)$ is the set of closed points of ${\cal Hom}(\Gamma,G)$ and there is a natural map
$$Hom(\Gamma,G)\to {\cal Hom}(\Gamma,G)$$
dual to
\begin{equation}\label{RandC}
R(\Gamma,G)\to R(\Gamma,G)/\sqrt{0}=\C[Hom(\Gamma,G)].
\end{equation}

By \cite[Thm 1.2]{KM}, $R(\Gamma,PSL(2,\C))$ contains non-zero nilpotent elements for some Artin groups $\Gamma.$ Furthermore, M. Kapovich proves that $R(\pi_1(M),SL(2,\C))$ and $R(\pi_1(M),PSL(2,\C))$ contain non-zero nilpotents for some $3$-dimensional manifolds $M,$ \cite{Ka1,Ka2}. See further comments in Sec. \ref{s_char_var_sch}.

\section{Spaces of irreducible representations}
\label{s_cr_reps}

\begin{proposition}\label{irred_open}
For every $\Gamma$ and every reductive group $G$ the set of irreducible $G$-representations
of $\Gamma$ is Zariski open in $Hom(\Gamma,G).$
\end{proposition}

\begin{proof}
The proposition follows from Corollary \ref{cr_irr} and from \cite[Prop 3.8]{Ne}.
Since the proof of this referenced result is non-elementary, we enclose a complete simple proof here:\\
(1) First, a simple proof for $G=GL(n,\C)$ and $SL(n,\C)$:
If $\rho: \Gamma\to G$ is irreducible then, by Shur's Lemma, the elements of $\rho(\Gamma)$ linearly span
the space of $n\times n$ matrices, $M(n,\C).$ Conversely, if $\rho(\Gamma)$ lies in a parabolic subgroup of $G$ then the elements of $\rho(\Gamma)$ do not span $M(n,\C).$

Enumerate all elements of $\Gamma$ in
a sequence $\gamma_1,\gamma_2,...$ Let $U_s$ be the space of all $\rho$'s such that
$\rho(\gamma_1),...,\rho(\gamma_s)$ span $M(n,\C).$ Since the space of all irreducible $\rho$'s
is the union of all $U_s$'s, it is enough to prove that each $U_s$ is open. This condition is equivalent
to an existence of a sequence $i_1,...,i_{n^2}$ such that the $n^2\times n^2$ matrix whose columns are $\rho(\gamma_{i_1}),...,\rho(\gamma_{i_{n^2}})$ considered as vectors in $M(n,\C)=\C^{n^2}$ has a non-zero determinant. This is a Zariski open condition.

(2) Here is a fairly elementary proof for all $G:$\\
The set of irreducible representations $\Gamma \to G$ is the complement of
$$\bigcup_P Hom(\Gamma,P)\subset Hom(\Gamma,G)$$
where the union on the left is over all proper parabolic subgroups of $G.$
By \cite[Thm 14.18]{Bo}, there are only finitely many parabolic subgroups of $G$ up to conjugation.
Therefore it is enough to prove that for a given $P$
$$X_P=\bigcup_{g\in G} Hom(\Gamma,gPg^{-1})\subset Hom(\Gamma,G)$$
is closed. $X_P$ is the union of closed sets parameterized by a complete variety $G/P.$
By Projective Extension Theorem, \cite[Ch 8 \S 5 Thm 6]{CLO}, such union is closed.
\end{proof}

The adjoint representation induces a map $Ad_*: Hom(\Gamma,G)\to Hom(\Gamma,GL(\g)).$
$\rho:\Gamma\to G$ is {\em $Ad$-irreducible} if $Ad\, \rho$ is irreducible.
Since the set of $Ad$-irreducible representations $\Gamma\to G$ is the $Ad_*$-preimage of the irreducible
representations in $Hom(\Gamma,GL(\g))$ we conclude with

\begin{corollary}\label{ad-irred_open}
The set of $Ad$-irreducible representations is Zariski open in $Hom(\Gamma,G)$.
\end{corollary}

\begin{proposition}\label{irred_exist}
Let $G$ be a reductive group.\\
(1) For a free group, $F_N,$ of rank $N\geq 2$,
the irreducible representations form a dense subset of $Hom(F_N,G)$ (in complex topology).\\
(2) For every surface group\footnote{A surface group is the fundamental group of a closed, orientable surface of genus $\geq 2$.} $\Gamma$, the irreducible
representations are dense in a non-empty set of irreducible components of $Hom(\Gamma,G).$
\end{proposition}

\begin{proof}
(1) Since $Hom(F_N,G)$ is an irreducible algebraic set and the set of irreducibles is Zariski open in it, it is enough to show that the set of irreducibles is non-empty.
Since every free group $F_N$ of rank $N\geq 2$ maps onto $F_2,$ it is enough to prove that statement for $F_2.$ The set of irreducible $G$-representations of $F_2$ is the complement of
$$\bigcup_P Hom(F_2,P)\subset Hom(F_2,G)=G\times G,$$
where the union of sets on the left is over all proper parabolic subgroups of $G.$
By \cite[Thm 14.18]{Bo}, there are only finitely many parabolic subgroups of $G$ up to conjugation.
Since for each $P$
$$\bigcup_{g\in G} Hom(F_2,gPg^{-1})$$
is the image of the $G$ action on $Hom(F_2,P)$ with stabilizer $P$,
its dimension is at most
$$2\cdot dim\, P+dim\, G-dim\, P< 2\cdot dim\, G=dim\, Hom(F_2,G)$$
Therefore, there exists an irreducible representation.\\
(3) Again, it is enough to prove that an irreducible representation exists. This follows from the fact that $\Gamma$ maps onto the free group of rank $2$.
\end{proof}

Note that Proposition \ref{irred_exist} does not apply to $F=$torus. For example, all representations of $\Z\times \Z$ to $GL(n,\C)$ are reducible for $n>1.$


%
\section{Stable and properly stable representations in the sense of GIT}
\label{s_git}
%

Let $O_\rho$ be the orbit of $\rho\in Hom(\Gamma,G)$ under the $G$ action on $Hom(\Gamma,G)$ by conjugation.
In the language of geometric invariant theory, $\rho$ is {\em poly-stable} if $O_\rho$ is closed.

\begin{theorem}\label{tfae-closed}
For any reductive algebraic group $G,$ $O_\rho\subset Hom(\Gamma,G)$ is closed if and only if $\rho$
is completely reducible.
\end{theorem}

\begin{proof} (The proof for $G=GL(n,\C)$, can be found in \cite[Thm 1.27]{LM})\\
$\Rightarrow$ We follow an argument of the proof of \cite[Thm 1.1]{JM}: Assume that $O_\rho$ is closed.
If $\rho$ is contained in a proper parabolic subgroup $P$ then by conjugating $\rho$ with a
one parameter group in the center of a Levi subgroup $L$ of $P$ one can obtain a representation $\rho'\in \overline O_\rho$ whose image lies in $L$. Since $O_\rho$ is closed, $\rho'=g^{-1}\rho g,$
for some $g\in P.$ Hence $\rho$ lies in the Levi subgroup $gLg^{-1}.$

$\Leftarrow$ Any finitely generated group $\Gamma$ is a quotient of a free group $F$ of finite rank. Denote the epimorphism $F\to \Gamma$ by $\pi.$ Since $Hom(\Gamma,G)$ is a closed subset of $Hom(F,G)$ and
$O_\rho=O_{\pi \rho}\cap Hom(\Gamma,G),$
it is enough to prove that \mbox{$O_{\pi \rho}\subset Hom(F,G)$} is closed.
This statement follows from \cite[Thm. 3.6]{Ri2}.
\end{proof}

According to the geometric invariant theory, a point $x$ of an affine set $X$ is {\em stable} with respect to a $G$ action on $X$ (and the trivial line bundle on $X$) if there is a Zariski open neighborhood of $x$ preserved by $G$ on which the $G$ action is closed, \cite{MFK, Do}.

\begin{cor}
(1) Every irreducible representation is a stable point of $Hom(\Gamma,G)$ under the $G$ action by conjugation.\\
(2) $\rho \in Hom(F_n,G)$ is stable if and only if $\rho$ is irreducible.
\end{cor}

\begin{proof}
(1) Follows from Proposition \ref{irred_open}.\\
(2) Every stable $\rho$ it is completely reducible by Theorem \ref{tfae-closed}. Every
completely reducible representation of a free group which is not irreducible can be deformed by
an arbitrarily small deformation to a representation which is not completely reducible.
\end{proof}

A point $x\in X$ is {\em properly stable} if it is stable and its stabilizer is finite.

\begin{corollary}\label{prop_stable} For every reductive group $G,$\\
(1) $\rho$ is a properly stable point of $Hom(\Gamma,G)$ under conjugation action of $G$ if and only if $\rho$ is irreducible and $C(G)$ is finite.\\
(2) $\rho$ is a properly stable point of $Hom(\Gamma,G)$ under conjugation action of $G/C(G)$ if and only if $\rho$ is irreducible.\\
\end{corollary}

\begin{proof} (1) $\Rightarrow:$ $\rho$ is completely reducible by Theorem \ref{tfae-closed}, and it is irreducible by Lemma \ref{centralizer}.\\
$\Leftarrow:$ by Theorem \ref{tfae-closed} and Propositions \ref{irred_open} and \ref{center}.\\
The same argument shows (2)
\end{proof}

Following \cite{JM}, we say that a representation $\rho$ is {\em good} if $O_\rho$ is closed and $S_\rho$ is the center of $G$. By Theorem \ref{tfae-closed} and Corollary \ref{cr_irr}, every good representation is irreducible. By Proposition \ref{center-ad}, every $Ad$-irreducible representation is good.

\begin{proposition}\label{g_open}
For every $\Gamma$ the space of good $G$-representations is open in the space of all
$G$-representations of $\Gamma$.
\end{proposition}

\begin{proof}
By \cite[Proposition 1.1]{JM}, the $G$ action on the space of all irreducible $G$-representations of $\Gamma$ is proper. The good representations, if they exist, form a set which is the the union of the principal orbits of that action. For every proper action, the union of principal orbits is an open subset, c.f. \cite[Thm. 1.5]{GO}.
\end{proof}

%
\section{Tangent Spaces}
\label{s_tangent}
%

Let $A$ be a commutative $\C$-algebra, let $m$ be a closed point of $Spec\, A$, i.e.
a maximal ideal $m\triangleleft A,$ and let $r_m$ be the projection $A\to A/m=\C.$
The tangent space to $Spec\, A$ at $m$ is the dual vector space to $m/m^2,$
$$T_m Spec\, A =(m/m^2)^*.$$

Here is an equivalent definition of the tangent space which will be useful later:
Let $\pi: \C[\ve]/(\ve^2) \to \C$ be the homomorphism sending $\ve$ to $0$ and
let ${\mathcal T}_m\, Spec\, A$ be the complex vector space of $\C$-algebra homomorphisms
$A\to \C[\ve]/(\ve^2)$ which descend to $r_m$ when composed with $\pi.$
Since every such homomorphism is of the form $r_m+\tau\ve,$ where $\tau: A\to \C$ is a derivation,
$${\mathcal T}_m\, Spec\, A=\{\tau: A\to \C: \tau(ab)=r_m(a)\tau(b)+r_m(b)\tau(a)\}.$$
A straightforward calculation shows that for every $v\in T_m\, Spec\, A,$
$$\lambda_v(a)=v(a-r_m(a))$$
is a derivation in ${\mathcal T}_m\, Spec A.$
A direct computation shows that
\begin{equation}\label{T-isom}
\lambda: T_m\, Spec\, A \to {\mathcal T}_m\, Spec\, A
\end{equation}
sending $v$ to $\lambda_v$ is an isomorphism of vector spaces, \cite[VI.1.3]{EH}.
From now on we will identify these two spaces and call them the {\em Zariski tangent space} to $Spec\, A$ at $m.$

The above discussion applies to $A={\cal Hom}(\Gamma,G)=Spec\, R(\Gamma,G).$ Each $\rho\in Hom(\Gamma,G)$ defines a projection
$r_\rho: R(\Gamma,G)\to \C$ and a closed point $m_\rho=Ker\, \rho$ in ${\cal Hom}(\Gamma,G).$
We will abbreviate $T_{m_\rho}\, {\cal Hom}(\Gamma,G)$ to $T_\rho\, {\cal Hom}(\Gamma,G).$
Each tangent vector
\mbox{$\tau\in T_\rho\, {\cal Hom}(\Gamma,G)$}
defines a group homomorphism
\begin{equation}\label{G-eps}
\Gamma\stackrel{\rho_U}{\longrightarrow} G(R(\Gamma,G)) \stackrel{G(r_\rho+\tau\ve)}{\longrightarrow} G(\C[\ve]/(\ve^2)).
\end{equation}

By abuse of notation, we denote by $\pi$ the extension of the homomorphism
$$\pi: \C[\ve]/(\ve^2)\to \C,\quad \pi(\ve)=0,$$
to the induced group homomorphism
$$\pi: G(\C[\ve]/(\ve^2))\to G(\C)=G.$$

\begin{proposition}\label{tangent}
Consider a closed embedding $G\subset GL(n,\C).$ (Such an embedding exists for
every affine algebraic group.)\\
(1) For every $g\in G(\C[\ve]/(\ve^2)),$
$$\sigma(g)=\frac{g\cdot \pi(g)^{-1}-I}{\ve}\in M_n(\C[\ve]/(\ve^2))$$ and
$\sigma(g)$ belongs to the Lie algebra $\g\subset M_n(\C)$ of $G$.\\
(2) For every $g_1,g_2\in G(\C[\ve]/(\ve^2))$,
\begin{equation}\label{cocycle}
\sigma(g_1g_2)=\sigma(g_1)+Ad\,\pi(g_1)\cdot \sigma(g_2),
\end{equation}
where $Ad: G\to GL(\g)$ is the adjoint representation.
\end{proposition}

\begin{proof}
(1) If $h\in G(\C[\ve]/(\ve^2))$ is such that $\pi(h)=I$ then $\frac{h-I}{\ve}$ belongs to the
Zariski tangent space to $G$ at the identity, that is the Lie algebra of $G.$
Now the statement follows by substitution $h=g\pi(g)^{-1}.$\\
(2) follows by a direct computation.
\end{proof}

For every $\tau\in T_\rho\, {\cal Hom}(\Gamma,G)$, $G(r_\rho+\tau\ve)\rho_U(\gamma)\in G(\C[\ve]/(\ve^2)),$
c.f. (\ref{G-eps}).
Therefore, by Proposition \ref{tangent}, we have a function $\sigma:\Gamma\to \g$
\begin{equation}\label{sigma}
\sigma(\gamma) = \frac{\left(G(r_\rho+\tau\ve)\rho_U(\gamma)\right)\left(G(r_\rho)\rho_U(\gamma)\right)^{-1}-I}
{\ve}
\end{equation}
satisfying (\ref{cocycle}), which is the cocycle condition
for the first cohomology group of $\Gamma$ with coefficients in $\g$ twisted by $Ad\, \rho.$
Hence, (\ref{sigma}) defines a linear map
$$\Psi_\rho: T_\rho\, {\cal Hom}(\Gamma,G)\to Z^1(\Gamma, Ad\,\rho)$$
sending $\tau$ to $\sigma.$

The adjoint action of the centralizer of $\rho(\Gamma)$, $S_\rho,$ on $\g$ induces a $S_\rho$-action on $Z^1(\Gamma, Ad\, \rho).$
Additionally, every $g\in S_\rho$ acts on $T_\rho\, {\cal Hom}(\Gamma,G)$ by sending $\tau\in
T_\rho\, {\cal Hom}(\Gamma,G)$ to $g\tau$ such that $$G(r_\rho+g\tau\ve)=gG(r_\rho+\tau\ve)g^{-1}.$$
The homomorphism $\Psi_\rho$ is a $S_\rho$-equivariant.

We are going to prove that $\Psi_\rho$ is an isomorphism by constructing its inverse.
An easy calculation shows that for every $\sigma\in Z^1(\Gamma,Ad\,\rho),$
$$\gamma\to (I+\sigma(\gamma)\ve)\cdot \rho(\gamma)$$
is a group homomorphism from $\Gamma$ to $G(\C[\ve]/(\ve^2))$ (c.f. \cite[Prop. 2.2]{LM} for $G=GL(n,\C)$). Therefore, $\sigma$ defines a homomorphism $\Phi_\rho(\sigma): R(\Gamma,G)\to \C[\ve]/(\ve^2)$
such that $\pi \Phi_\rho(\sigma)=r_\rho.$ Hence, $\Phi_\rho(\sigma)\in T_\rho\, {\cal Hom}(\Gamma,G).$
In other words, we have defined a linear map
$$\Phi_\rho: Z^1(\Gamma, Ad\, \rho) \to T_\rho\, {\cal Hom}(\Gamma,G).$$
A straightforward computation shows (c.f. \cite[Lemma 2.2]{JM}, \cite[Prop 2.2]{LM} for $G=GL(n,\C)$ and \cite[Prop 1.2]{Be} for $G=PSL(2,\C)$):

\begin{theorem}\label{Z^1-thm}
$\Psi_\rho$ and $\Phi_\rho$ are inverses of each other, and therefore, they are $S_\rho$-equivariant isomorphisms between
$Z^1(\Gamma,Ad\, \rho)$ and $T_\rho\, {\cal Hom}(\Gamma,G).$
\end{theorem}

%
\section{Smooth and Reduced Representations}
\label{s_reduced_rep}
%

A closed point $x$ of an algebraic scheme $X$ is {\em reduced} if the local ring $O_{X,x}$ has no non-zero nilpotent elements.
By this definition, reduced points form a Zariski open subset of $X.$

Recall that a point $x$ of an affine algebraic set or of an algebraic scheme $X$ is simple if $x$ belongs to a unique irreducible component of $X$ and the dimension of that component coincides with $dim\, T_x\, X.$
Every simple point of an algebraic scheme is reduced.

We say that $\rho:\Gamma\to G$ is {\em smooth} (respectively: {\em scheme smooth}), if $\rho$ is a simple point of $Hom(\Gamma,G)$
(respectively: of ${\cal Hom}(\Gamma,G)$).
$\rho$ is {\em reduced} if $\rho$ is a reduced point of ${\cal Hom}(\Gamma,G).$
Note that $\rho$ is scheme smooth iff it is reduced and smooth.

\begin{remark}\label{reduced_open}
(1) The set of all smooth representations,\\
(2) the set of all reduced representations,\\
(3) the set of all scheme smooth representations,\\
are Zariski open subsets $Hom(\Gamma,G).$
\end{remark}

\begin{proof} (1) follows from the fact that the set of simple points of Zariski open.\\
(2) $Hom^r(\Gamma,G)$ is a preimage of the Zariski open set of reduced points in ${\cal Hom}(\Gamma,G)$
under the map $Hom(\Gamma,G)\to {\cal Hom}(\Gamma,G).$\\
(3) This set is the intersection of the first two.
\end{proof}

All $G$-representations of a free group are scheme smooth, since
$R(F_n,G)$ is the coordinate ring of the $n$-th cartesian power of $G,$
which is a non-singular algebraic set.

\begin{proposition}\label{surface-reduced}
For every reductive group $G$ and every surface group $\Gamma$ all irreducible representations $\rho:\Gamma\to G$ are scheme smooth and, hence, reduced.
\end{proposition}

\begin{proof}
By Proposition \ref{center}, the centralizer of $\rho(\Gamma)$ is a finite extension of the center of $G$.
Hence, by \cite[Prop. 1.2]{Go1} and by Theorem \ref{Z^1-thm},
\begin{equation}\label{ineq}
dim\, T_\rho\, {\cal Hom}(\Gamma,G)=dim\, Z^1(\Gamma, Ad\, \rho)=(2g-1) dim\, G+ dim\, C(G).
\end{equation}

(1) Assume first that $G$ is semi-simple. Then $dim\, C(G)=0.$
Since $\Gamma$ has a presentation with $2g$ generators and one relation,
$$dim\, C\geq (2g-1) dim\, G$$
for all irreducible components $C\subset Hom(\Gamma,G).$
Therefore,
\begin{equation}\label{dimC}
dim\, C\geq dim\, T_\rho\, {\cal Hom}(\Gamma,G).
\end{equation}
Since (\ref{dimC}) must be an equality, $\rho$ is scheme smooth.

(2) For an arbitrary reductive group $G$ consider epimorphism (\ref{nu}),
$$\nu: C^0(G)\times [G,G] \to G.$$
Since it has a finite kernel, the induced map
$$Hom(\Gamma,C^0(G))\times Hom(\Gamma,[G,G])\to Hom(\Gamma,G)$$
is finite.
Since $Hom(\Gamma,C^0(G))=\left(C^0(G)\right)^{2g}$ and, by (1), $Hom(\Gamma,[G,G])$ is composed of
irreducible components of dimension at least $(2g-1)\dim\, [G,G],$ the set
$Hom(\Gamma,G)$ is composed of irreducible components of dimension at least
$$2g\cdot dim\, C^0(G)+(2g-1)\dim\, [G,G]=(2g-1) dim\, G+ dim\, C^0(G).$$
Therefore, $$dim\, C\geq dim\, T_\rho\, {\cal Hom}(\Gamma,G).$$
Now the argument goes exactly as in (1).
\end{proof}

By \cite[II \S2 Thm 6]{Shf}, we obtain:

\begin{corollary} For every reductive group $G$, every irreducible representation of every surface group $\Gamma$ belongs to a unique irreducible component of $Hom(\Gamma,G).$
\end{corollary}

\begin{proposition}\label{surf-R-reduced}
If $G$ is a reductive group and $\Gamma$ is a surface group such that $Hom(\Gamma,G)$ is irreducible, then $R(\Gamma,G)$ has no non-zero nilpotents. Consequently, map (\ref{RandC}), $R(\Gamma,G)\to \C[Hom(\Gamma,G)],$ is an isomorphism.
\end{proposition}

\begin{proof} Recall from Sec. \ref{s_tangent} that for every $\rho:\Gamma\to G$ there is a corresponding maximal ideal $m_\rho \triangleleft R(\Gamma,G).$
If we assume that $R(\Gamma,G)$ has a non-zero nilpotent, then that nilpotent projects onto a non-zero nilpotent in the localization $R(\Gamma,G)_{m_\rho}$ for an open set of
$\rho$'s in $Hom(\Gamma,G).$ These $\rho$'s are not reduced.

Since $Hom(\Gamma,G)$ is irreducible, by Proposition \ref{irred_exist}(2), the set of irreducible representations is dense in $Hom(\Gamma,G)$. Therefore, $\Gamma$
has an irreducible but not reduced $G$-representation, contradicting Proposition \ref{surface-reduced}.
\end{proof}

A version of this theorem for $F$ torus appears in \cite{Si-ab}.

\begin{corollary}\label{no_nilp_gl}
For $G=GL(n,\C)$ and $SL(n,\C)$ and every surface group $\Gamma$ the ring $R(\Gamma,G)$ is an integral domain, i.e. it has no proper zero divisors.
\end{corollary}

\begin{proof}
By Theorems 1 and 3 in \cite{RBC}, $Hom(\Gamma,GL(n,\C))$ and $Hom(\Gamma,SL(n,\C))$ are irreducible. Therefore, the statement follows by Proposition \ref{surf-R-reduced}.
\end{proof}

This corollary has some important consequences to the theory of skein modules, c.f.
Corollary \ref{skein}, and it is needed to complete Goldman's construction of symplectic forms on character varieties of surfaces, c.f. Section \ref{s_char_var_surf}.

The following remains open:

\begin{question}
Does $R(\Gamma,G)$ have non-zero nilpotents for some algebraic group $G$ and a surface group $\Gamma$?
\end{question}

Finally, we would like to remark that Proposition \ref{surface-reduced} does not hold for non-surface groups.
In fact, there appears to be no easy characterization of simple points of $Hom(\Gamma,G)$ in general.


\begin{example}
Let $\rho_1:\Z^2\to SL(2,\C)$ be the trivial representation and let $\rho_2:\Gamma \to SL(2,\C)$ be an irreducible representation.
These representations define a representation $\rho_1*\rho_2: \Z^2*\Gamma\to SL(2,\C)$ which is irreducible.
On the other hand, the Jacobian matrix $(\p r_i/\p x_j)_{i=1,...,5,j=1,...,8}$ of the five relations
in Example \ref{Z^2} has rank $2$ at $\rho_1=(x_1,...,x_8)=(1,0,0,1,1,0,0,1)$.
Therefore, $Hom(\Z^2,SL(2,\C))$ is singular at $\rho_1$ and, by Remark (\ref{free_prod}), $\rho_1*\rho_2$ is a singular point of $Hom(\Z^2*\Gamma,G).$
\end{example}

%
\section{Orbits}
\label{s_orbits}
%

As before, let $O_\rho$ be the orbit of $\rho$ in $Hom(\Gamma,G)$ under the $G$ action by conjugation.
Since $O_\rho$ is homogeneous and (as every algebraic set) it has a simple point, all its points are simple, i.e.
$O_\rho$ is a non-singular algebraic set.

The following theorem generalizes \cite[Lemma 2.2]{JM} and \cite[Cor 2.4]{LM}.

\begin{thm}\label{orbit}
For every $\rho$ the inclusion
$$T_\rho\, O_\rho\subset T_\rho\, {\cal Hom}(\Gamma,G)$$
corresponds to $$B^1(\Gamma,Ad\, \rho) \subset Z^1(\Gamma,Ad\, \rho)$$
under the isomorphism $\Psi_\rho.$
\end{thm}

\begin{proof}
Since $O_\rho$ is the image of the map $f_\rho: G\to Hom(\Gamma,G),$ $f_\rho(g)= g\rho g^{-1}$,
$O_\rho$ is the left quotient of $G$ by the stabilizer of $\rho,$ $S_\rho,$ c.f. \cite[II.6.1]{Bo}.
Furthermore, the quotient map $G\to G/S_\rho$ induces an epimorphism $T_{[\rho]}\, G\to T_\rho\, (G/S_\rho).$
(Proof: $S_\rho$ acts on $G$ by left translation, with the trivial stabilizer. By Luna's \'Etale Slice Theorem, \cite{Lu},
$G\to G//S_\rho$ is \'etale, inducing an isomorphism of tangent spaces. Note that $G//S_\rho=G/S_\rho$ since all orbits are closed.) Therefore every tangent vector $v$ to $O_\rho$ at $\rho$ is of the form $d(f_\rho)(I+A\ve).$
By the discussion of Section \ref{s_tangent}, every vector $v\in T_\rho\, {\cal Hom}(\Gamma,G)$ corresponds to a group homomorphism
(\ref{G-eps}). By this correspondence, $v$ corresponds to
\begin{equation}\label{tangent_rep2}
\gamma\to(I+A\ve)\rho(\gamma)(I+A\ve)^{-1}.
\end{equation}
Since
$$(I+A\ve)^{-1}=I-A\ve\quad \text{mod\ } \ve^2,$$
(\ref{tangent_rep2}) can be written as
$$\gamma\to (I+A\ve)\rho(\gamma)(I-A\ve)=(I+(A-\rho(\gamma)A\rho(\gamma)^{-1})\ve)\rho(\gamma)\quad {mod}\ \ve^2.$$
The right side is $(I+\tau\ve)\rho(\gamma),$ where
\begin{equation}\label{tau}
\tau=A-\rho(\gamma)A\rho(\gamma)^{-1}\in B^1(\Gamma,Ad\, \rho)\subset Z^1(\Gamma,Ad\, \rho).
\end{equation}
Therefore, $\Psi_\rho$ of Theorem \ref{Z^1-thm} maps $T_\rho\, O_\rho$ to $B^1(\Gamma,Ad\, \rho).$
Since $\Psi_\rho$ is an embedding, in order to show $\Psi_\rho$ is an isomorphism from $T_\rho\, O_\rho$ to $B^1(\Gamma,Ad\, \rho)$, it is enough to prove that
\begin{equation}\label{TOB}
dim\, T_\rho\, O_\rho \geq dim\, B^1(\Gamma,Ad\, \rho).
\end{equation}
Since $Ad\rho(\gamma):\g\to \g$ is the differential of the map
$$G\to G,\quad g\to \rho(\gamma)g \rho(\gamma)^{-1}$$
and that map is constant on $S_\rho,$ $Ad\rho(\Gamma)$ acts trivially on the Lie algebra $L(S_\rho)\subset \g.$
Therefore (\ref{tau}) vanishes for $A\in L(S_\rho).$ Consequently,
$$dim\, B^1(\Gamma,Ad\, \rho)\leq \dim\, \g-\dim\, S_\rho.$$
Since the left side of (\ref{TOB}) is $dim\, G-dim\, S_\rho,$ inequality (\ref{TOB}) follows.
\end{proof}

\begin{remark}\label{B^1-equiv}
The action of $S_\rho$ on $T_\rho\, {\cal Hom}(\Gamma,G)$ and on $Z^1(\Gamma,Ad\, \rho)$ preserves $T_\rho\, O_\rho$ and $B^1(\Gamma,Ad\, \rho).$
\end{remark}


Theorems \ref{Z^1-thm} and \ref{orbit} imply Weil's infinitesimal rigidity theorem, c.f. \cite{Wei}.

\begin{corollary}\label{local_rig}
If $H^1(\Gamma,Ad\, \rho)=0$  then\\
(1) $\rho:\Gamma\to G$ is scheme smooth (and, hence, reduced)\\
(2) $O_\rho$ contains an open neighborhood of $\rho$ (in complex topology) in $Hom(\Gamma,G)$.
\end{corollary}

\begin{proof}
(1) Since $O_\rho$ is smooth, we have
$$dim\, B^1(\Gamma,Ad\, \rho)=dim\, O_\rho\leq dim\, Hom(\Gamma,G)\leq
dim\, T_\rho Hom(\Gamma,G)$$
$$\leq dim\, T_\rho {\cal Hom}(\Gamma,G)= dim\, Z^1(\Gamma,Ad\,\rho),$$
by Theorems \ref{Z^1-thm} and \ref{orbit}.
Since all relations above are equalities, $\rho$ is scheme smooth.\\
(2) follows from the fact that $Hom(\Gamma,G)$ is smooth at $\rho$
and has the same dimension as $O_\rho$.
\end{proof}

Since $H^1(\Gamma,Ad\,\rho)=0$, for all finite groups $\Gamma$ and all $\rho$, c.f. \cite[Thm, 6.5.8]{Wb}, we also have:

\begin{corollary}\label{finite_reduced}
All $G$-representations of every finite group are reduced.
\end{corollary}

%
\section{Character Varieties}
\label{s_char_var}
%

The categorical quotient of $Hom(\Gamma, G)$ by the $G$ action by conjugation,
$$X_G(\Gamma)= Hom(\Gamma, G)//G,$$
is called the {\em $G$-character variety} of $\Gamma.$
By definition, it is an affine algebraic set together with the map
$Hom(\Gamma,G)\to Hom(\Gamma, G)//G$ which is constant on all $G$-orbits, with the property that every morphism from $Hom(\Gamma,G)$ into an affine algebraic set $Y$ which is constant on all $G$-orbits factors through $Hom(\Gamma,G)\to Hom(\Gamma,G)//G,$ c.f. \cite{Do,Fo,MFK}.
If $G$ is reductive then the categorical quotient exists.
The reason for considering the categorical quotient rather than the set theory quotient is that the quotient topology on $Hom(\Gamma, G)/G$ is often not a Zariski topology of any algebraic set. In particular, it often contains points which are not closed.
Character varieties are often reducible, despite being called ``varieties".

Every equivalence class in $X_G(\Gamma)$ contains a unique closed orbit. Therefore, by
Proposition \ref{tfae-closed}, each element of the $G$-character variety of $\Gamma$ is represented by
a unique completely reducible representation.

\begin{example}
Let $\T$ be a maximal torus of $G.$ The map $\T\to Hom(\Z,G)$ sending $g$ to the $G$-representation
$1\to g$ of $\Z$ factors through an isomorphism $$\T/W\to X_G(\Z),$$
where $W$ is the Weyl group of $G,$ c.f. \cite[6.4]{St}.
\end{example}


\begin{example}
(1) The $SL(2,\C)$-character variety of the free group on two generators is isomorphic to $\C^3.$\\
(2) $SL(3,\C)$-character variety of the free group on two generators is a hypersurface in $9$-dimensional affine space, \cite[Thm 8]{La1}, \cite{Si1}.
\end{example}


Algebraic properties of character varieties are discussed in \cite{AP,BH,BK1,BK2,BK3,BKCh,GM,Ho1,Ho2,La1,La2,La3,La4,La5,LP,LM,Na, PBK,RBK,RBC,Si1,Si-g,Si-ab,Wh}.
Character varieties appear in many ways in mathematics and physics.
Of particular importance are character varieties of surface groups, discussed in Section \ref{s_char_var_surf}. They appear as moduli spaces of hyperbolic, projective, and other geometric structures on surfaces, c.f. \cite{BIW,Go1,Go3,Go9,GM1,GM2,GW,JM,KM,Li,Sa,Wa},
as well as and the moduli spaces of flat connections, holomorphic bundles, and of Higgs bundles, \cite{AB,AMW,CHM,Da,DDW,DWWW,DWW,GGM,FL1,FL2,GGM,Hi1,Hi2,HLR,HT,Je2,Je3,JK,Ki,NS,Me,MW,Ol,Rac,
Sim1,Sim2,Th1,Th2,Th3,Wi2,Wi3,Za}. These connections inspired an investigation of topology (and, more specifically, cohomology) of character varieties, in many of the papers cited here.

In a related fashion, character varieties of surfaces appear also as examples of symplectic reduction (and, furthermore, K\"ahler and hyper-K\"ahler reductions) as well as in the context of Hamiltonian actions, \cite{AB,Au,Go-symp,Je1}.
The action of the mapping class groups on character varieties is discussed in
\cite{Go6,Go7,Go8,PX,SS,We,Wi}.

In mathematical physics, character varieties appear in the context of Yang-Mills and Chern-Simons quantum field theories, \cite{At,AB,Ba,Fr,Gu,JW,KK,We1,We2,We3,Wi1,Wi2,Wi3}
as well as in the related skein theory of quantum invariants of $3$-manifolds, \cite{Bu,CM,FG,FGL,Ga,Ge,Le2, PS1,PS2,Si2,Si3}.

In low-dimensional topology, character varieties appear in the Bass-Culler-Shalen theory, in the context of A-polynomial and in other related areas, \cite{Be,BB,BLZ,BZ1,BZ2,BZ3,Du,DG,CCGLS,CL,CS,HS,MS,Mo,Shl} as well as
in the context of Casson-Walker-Lescop invariant, \cite{AM,BC1,BC2,BHe,BN,Cu}.
Varieties of representations of $3$-manifold groups have been studied for a long time also in relation to the Alexander polynomial, hyperbolic geometry, and for other independent reasons, \cite{BF,Ril1,HLM1,HLM2,HP1,HP2,Le1,LR1,LR2,Po,Ti}.
Finally, character varieties appear in the discussion of local rigidity of discrete subgroups of Lie groups, c.f. Corollary \ref{local_rig}, \cite{Wei,Ra,Rag}.

These are just sample references to the above topics, as there are hundreds of papers
devoted to every one of them.

Denote the set of irreducible representations in $Hom(\Gamma,G)$ by $Hom^i(\Gamma,G).$
The $G$ action by conjugation preserves this set. Since all orbits in $Hom^i(\Gamma,G)$ are closed (c.f. Proposition \ref{tfae-closed}) and each equivalence class in a
categorical quotient contains a unique closed orbit, the categorical quotient,
$Hom^i(\Gamma,G)//G,$ is the set-theoretic quotient as well.
Denote it by $X^i_G(\Gamma).$

\begin{proposition}\label{orbifold}
Let $G$ be a reductive group.\\
(1) For the free group $F_n,$ $X^i_G(F_n)$ is a complex orbifold of complex dimension $$(n-1)dim\, G+dim\, C(G).$$
(2) For every closed orientable surface $S_g$ of genus $g>1,$ $X^i_G(\pi_1(S_g))$ is a complex orbifold of complex dimension $$(2g-2)dim\, G+2dim\, C(G).$$
(3) If $G$ is CI (e.g. $G=GL(n,\C), SL(n,\C)$) then $X_G^i(F_n)$ and $X_G^i(\pi_1(S_g))$ are manifolds for all $n,g>1$. (See also \cite{FL2}.)
\end{proposition}

\begin{proof}
(1) Since $Hom^i(F_n,G)$ is an open subset of $G^n,$ it is smooth.
By \cite[Prop 1.1]{JM}, the $G/C(G)$ action on $Hom^i(F_n,G)$ is properly discontinuous. By Proposition \ref{center}, $G/C(G)$ acts on $Hom^i(F_n,G)$ with finite centralizers.
The quotient of a manifold by a properly discontinues group action with finite stabilizers is an orbifold.\\
(2) By Proposition \ref{surface-reduced}, all irreducible representations of $\pi_1(S_g)$ are scheme smooth. Therefore, by Theorem \ref{Z^1-thm} and by \cite[Prop. 1.2]{Go-symp}, $Hom^i(\pi_1(S_g))$ is a complex manifold of dimension
$(2g-1)dim\, G+dim\, C(G).$  Now the proof follows as in (1).\\
(3) By definition of a CI group, the $G/C(G)$ action on $X_G^i(F_n)$ is free. By \cite[Prop 1.1]{JM}, this action is also properly discontinuous. The quotient of a manifold by a free properly discontinues group action is a manifold, \cite[Proposition 3.5.7]{Th}.
\end{proof}

Recall that a representation $\rho: \Gamma\to G$ is good if and only if it is irreducible and the centralizer of $\rho(\Gamma)$ coincides with the center of $G.$ By Proposition \ref{g_open}, the space of good representations, $Hom^g(\Gamma,G),$ is an open subset of the irreducible ones. Since $G/C(G)$ acts freely and properly discontinuously on $Hom^g(\Gamma,G)$, we have:

\begin{corollary}\label{good-smooth}
For every reductive group $G$ and every surface group or a free group $\Gamma$,
$X^g_G(\Gamma)=Hom^g(\Gamma,G)/G$ is an open subset of $X^i_G(\Gamma)$ and a smooth complex manifold.
\end{corollary}

For a topological space $Y,$ we will abbreviate $X_G(\pi_1(Y))$ by $X_G(Y).$

\begin{proposition}\label{3mfld}
(1) If $M$ is a compact $3$-manifold with a connected boundary of genus $g$ then
$$dim\, X_G(M)\geq dim\, G\cdot (g-1)+dim\, C(G).$$
(2) For a given non-abelian reductive group $G$ and a positive integer $g$ there is no upper bound on $dim\, X_G(M)$ over compact $3$-manifolds $M$ with connected boundary of fixed genus $g.$
\end{proposition}

\begin{proof}
(1) If $M$ is a compact manifold with connected boundary of genus $g$ then
$\pi_1(M)$ has a presentation with $n$ generators and $p$ relations such that
$$1-n+p=\chi(M)=1-g.$$
Hence $dim\, Hom(\pi_1(M),G)\geq dim\, G\cdot (n-p)=dim\, G\cdot g.$
Since $X_G(M)$ is the quotient of $Hom(\pi_1(M),G)$ by $G/C(G),$ the statement follows.

(2) It is enough to prove that there is no upper bound on $dim\, Hom(\pi_1(M),G),$ over compact manifolds $M$ with connected boundary of fixed genus $g.$
Since every non-abelian reductive group contains either $SL(2,\C)$ or $PSL(2,\C),$
$$dim\, Hom(\pi_1(M),G)\geq min(dim\, Hom(\pi_1(M),SL(2,\C)), dim\, Hom(\pi_1(M),PSL(2,\C)).$$
Since the quotient map $SL(2,\C)\to PSL(2,\C)$ induces a finite map
$$Hom(\pi_1(M),SL(2,\C))\to Hom(\pi_1(M),PSL(2,\C)),$$
$$dim\, Hom(\pi_1(M),SL(2,\C))\leq dim\, Hom(\pi_1(M),PSL(2,\C)).$$
(The inequality stems from the fact that this map does not have to be onto.)
Therefore, it is enough to prove that there is no upper bound on $dim\, Hom(\pi_1(M),SL(2,\C)).$

Let $K_n$ be the connected sum of $n$ copies of a knot $K.$ Cooper and Long, \cite{CL}, proved that $dim\, Hom(\pi_1(S^3\setminus K_n), SL(2,\C))\geq n+3.$ (Although their argument is made for hyperbolic knots $K$ only, it generalizes to all knots by the result of \cite{KrM}, c.f. \cite{DG}.) Let $K_{n,g}$ be a graph obtained by connecting $g$ unlinked copies of $K_n$ in $S^3$ by $g-1$ tunnels and let $M_{n,g}$ be the complement of an open tabular neighborhood of $K_{n,g}$ in $S^3.$ Then $\pi_1(M_{n,g})$ is the free product of $g$ copies of $\pi_1(S^3\setminus K_n)$ and
$$dim\, Hom(\pi_1(M_{n,g}),SL(2,\C))\geq g\cdot (n+3),$$
by Remark \ref{free_prod}.
Since $\p M_{n,g}$ is a surface of genus $g$, the statement follows.
\end{proof}


\section{Character Varieties as Algebraic Schemes}
\label{s_char_var_sch}

The algebraic scheme ${\cal X}_G(\Gamma)={\cal Hom}(\Gamma,G)//G$ is often also called the $G$-character variety of $\Gamma.$
By the definition of categorical quotient, ${\cal X}_G(\Gamma)=Spec\, R(\Gamma,G)^G.$
The epimorphism $R(\Gamma,G)\to \C[Hom(\Gamma,G)]$
induces the map
\begin{equation}\label{nil-q}
R(\Gamma,G)^G\to \C[Hom(\Gamma,G)]^G.
\end{equation}
which is an epimorphism by the properties
of the Reynolds operator. Hence, we have an embedding of schemes
\begin{equation}\label{spec_X}
i: X_G(\Gamma)\hookrightarrow {\cal X}_G(\Gamma).
\end{equation}
Since the kernel of (\ref{nil-q}) is the nil-radical of $R(\Gamma,G)^G,$
(\ref{spec_X}) is a bijection between the points of $X_G(\Gamma)$
and the closed points of ${\cal X}_G(\Gamma).$

In \cite{Si1}, we have described ${\cal X}_{SL(n,\C)}(\Gamma)$ as a space of $n$-valent
graphs reminiscent of Feynman diagrams in an arbitrary path connected topological space $X$ with $\pi_1(X)=\Gamma.$

By Corollary \ref{no_nilp_gl}, $R(\Gamma,G)^G$ has no zero divisors for surface groups
and $G=GL(n,\C), SL(n,\C)$. Hence, by \cite{RBC}, ${\cal X}_G(\Gamma)$ is reduced and irreducible for such groups. For $F$ torus and $G=SL(2,\C)$, $R(\pi_1(F),G)^G$ is reduced as well by \cite[Thm 3.3]{PS1}.
These facts have an important consequence for the theory of skein modules.

\begin{corollary}\label{skein} For every closed orientable surface $F$,
the map $\phi$ of \cite[Thm 7.1]{PS2} is an isomorphism between the skein
algebra of $\pi_1(F)$ and $\C[X_{SL(2,\C)}(\pi_1(F))].$
Consequently, $\phi$ composed with the isomorphism of \cite[Thm 2.8]{PS2} is an isomorphism from the skein algebra of a surface,
$S_{2,\infty}(F,\C,-1),$ to $\C[X_{SL(2,\C)}(\pi_1(F))].$
\end{corollary}

Although this result was announced in \cite[Thm 7.3]{PS2}, its proof required
\cite[Thm 4.7]{PS2} whose proof was not provided.
An alternative proof of the above statement was provided independently by L. Charles and J. March\'e in \cite{CM}.

${\cal X}_G(\Gamma)$ is not always reduced. Kapovich and Millson proved that for every affine (possibly unreduced) variety $X$ over $\Q$
there is an Artin group $\Gamma$ such that a Zariski open subset of ${\cal X}_{PSL(2,\C)}(\Gamma)$ is isomorphic to a Zariski open subset of $X,$ \cite{KM}.
Additionally, for every $x\in X$ there is a representation $\rho$ of an Artin group $\Gamma$ into $PSL(2,\C)$ such that the analytic germ of ${\cal X}_{PSL(2,\C)}(\Gamma)$
at $[\rho]$ coincides with the analytic germ of $X$ at $x$, \cite{KM}.

Kapovich proved that the same is true for $3$-manifold groups. That is, for every $x\in X$ as above
there is a closed $3$-manifold $M$ and a $PSL(2,\C)$-representation $\rho$ of $\Gamma=\pi_1(M)$ such that the analytic germ of ${\cal X}_{PSL(2,\C)}(\Gamma)$ at $[\rho]$ coincides with the analytic germ of $X$ at $x$, \cite{Ka1,Ka2}.
In particular ${\cal X}_{PSL(2,\C)}(\Gamma)$ contains non-zero nilpotent elements for some Artin groups and some $3$-manifold groups $\Gamma.$

%
\section{Tangent spaces to character varieties}
\label{s_tan_char}
%

For every $\rho$ the map (\ref{spec_X}) induces an embedding
\begin{equation}\label{embedding}
T_{[\rho]}\, X_G(\Gamma)\hookrightarrow T_{[\rho]}\, {\cal X}_G(\Gamma).
\end{equation}
It is an isomorphism if $\rho$ is reduced. We are going to give a cohomological description of these tangent spaces.

For every $\rho:\Gamma\to G,$ the $S_\rho$ action on $Z^1(\Gamma,Ad\, \rho)$ descends to
an action on $H^1(\Gamma,Ad\, \rho).$

\begin{theorem}\label{tangent_X} (1) If $\rho:\Gamma\to G$ is completely reducible then the isomorphism $\Phi_\rho$
of Theorem \ref{Z^1-thm} combined with the natural projection ${\cal Hom(\Gamma,G)}\to {\cal X}_G(\Gamma)$ induces a natural
linear map $$\phi: H^1(\Gamma,Ad\, \rho)\to T_{[\rho]}\, {\cal X}_G(\Gamma).$$
(2) If $\rho$ is good then $\phi$ is an isomorphism.\\
(3) If $\rho$ is completely reducible and scheme smooth (but not necessarily good) then
$$dim\, T_0\, \left( H^1(\Gamma,Ad\, \rho)//S_\rho \right)=dim\, T_{[\rho]}\, {\cal X}_G(\Gamma)=dim\, T_{[\rho]}\, X_G(\Gamma).$$
\end{theorem}

The naturality of the morphism $\phi$ means that for every $\alpha:\Gamma'\to \Gamma$ and $\beta: G\to G'$ such that $\rho:\Gamma\to G$ and
$\beta\rho\alpha: \Gamma'\to G'$ are completely reducible, the following diagram commutes
$$\begin{array}{ccc}
H^1(\Gamma,Ad\, \rho) & \stackrel{\phi}{\longrightarrow} &  T_{[\rho]} {\cal X}_G(\Gamma)\\
\downarrow & & \downarrow\\
H^1(\Gamma',Ad\, \beta\rho\alpha)  & \stackrel{\phi}{\longrightarrow} & T_{[\beta\rho\alpha]} {\cal X}_{G'}(\Gamma'),
\end{array}$$
where the vertical homomorphisms are induced by $\alpha$ and $\beta.$

We discuss the existence of a natural isomorphism in part (3) of the above theorem at the end of this section.

Versions of the above theorem for $G=PSL(2,\C)$ appear in \cite[Prop 3.5]{Po} and \cite[Prop. 5.2]{HP2}.
(Note that every irreducible $PSL(2,\C)$-representation is $Ad$-irreducible and, hence, good, by Proposition \ref{center-ad}.) A related discussion (for abelian representations) appears in \cite{Be}.

Although statements similar to Theorem \ref{tangent_X}(2) appear in literature for general $G$, there is a lot of confusion about the necessary assumptions and there is a lack of discussion of the naturality of the isomorphism. Furthermore, all proofs known to us are incomplete.

\noindent{\it Proof of Theorem \ref{tangent_X}:}
The proof of the above theorem is based on Luna's \'Etale Slice Theorem, \cite{Lu}, c.f. \cite{MFK},\cite[Thm 6.1]{PV},
which applies to closed orbits of group actions on affine (not necessarily reduced) schemes.\\
(1) Since the orbit of $\rho$ is closed by Theorem \ref{tfae-closed}, by Luna's Slice Theorem there
exists a subscheme $S\subset {\cal Hom}(\Gamma,G),$ called an \'etale slice,
such that
\begin{enumerate}
\item[(i)] $\rho\in S$
\item[(ii)] $S$ is preserved by the $S_\rho$ action.
\item[(iii)] if we denote the categorical quotient of $G\times S$ by the diagonal $S_\rho$-action
(which is right on $G$ and left on $S$) by $G\times_{S_\rho} S$ then
the natural map
$G\times_{S_\rho} S\to {\cal Hom}(\Gamma,G)$
sending $(g,s)$ to $gs$ is \'etale.
\item[(iv)]
$S//S_\rho\to  {\cal X}_G(\Gamma)$
is \'etale.
\end{enumerate}
In particular, it means that the differentials of the above maps
\begin{equation}\label{Tetale1}
T_{(e,\rho)}\, \big(G\times_{S_\rho} S\big)\to T_\rho\, {\cal Hom}(\Gamma,G)
\end{equation}
and
\begin{equation}\label{Tetale2}
T_{[\rho]}\, S//S_\rho\to T_{[\rho]} {\cal X}_G(\Gamma)
\end{equation}
are isomorphisms.

The natural projection,
$${\cal Hom}(\Gamma,G)\to {\cal X}_G(\Gamma)$$
induces
\begin{equation}\label{Tpro}
T_\rho\, {\cal Hom}(\Gamma,G)\to T_{[\rho]}\, {\cal X}_G(\Gamma).
\end{equation}
which can be written as
\begin{equation}\label{prephi}
T_\rho\, {\cal Hom}(\Gamma,G)\stackrel{\simeq}{\longrightarrow} T_{(e,\rho)}\, \big(G\times_{S_\rho} S\big)\to
T_{[\rho]} (S//S_\rho)\stackrel{\simeq}{\longrightarrow} T_{[\rho]}\, {\cal X}_G(\Gamma),
\end{equation}
where the first map is the inverse of (\ref{Tetale1}), the second is the differential of
\begin{equation}\label{gs}
G\times_{S_\rho} S\to S//S_\rho,
\end{equation}
and the third one is (\ref{Tetale2}).

Since the map $G\to G\times_{S_\rho} S$ composed with (\ref{gs}) is constant,
the image of $T_e\, G$ in $T_{(e,\rho)} G\times_{S_\rho} S$ belongs to the kernel of $T_\rho\, {\cal Hom}(\Gamma,G)\to T_\rho (S//S_\rho)$. Since this image corresponds to $T_\rho\, O_\rho\subset T_\rho\, {\cal Hom}(\Gamma,G),$
(\ref{prephi}) factors through
\begin{equation}\label{prephi2}
T_\rho\, {\cal Hom}(\Gamma,G)/T_\rho\, O_\rho \to T_\rho (S//S_\rho)\stackrel{\simeq}{\longrightarrow} T_\rho\, {\cal X}_G(\Gamma).
\end{equation}
Now the statement follows from Theorems \ref{Z^1-thm} and \ref{orbit}.

(2) If $\rho$ is good then $S_\rho=C(G)$ acts trivially on $S$ and (\ref{gs}) can be written as the projection $$G/C(G)\times S\to S.$$
Its differential is the projection $T_e\, G/C(G) \times T_\rho\, S\to T_\rho\, S$
whose kernel is $T_\rho\, O_\rho.$
Therefore, the left morphism of (\ref{prephi2}) is an isomorphism.

(3) Since the stabilizer of the $S_\rho$-action on $G\times S$ at $(e,\rho)$ is trivial, by Luna's Slice Theorem
we have an \'etale morphism $(G\times S)/\{e\}\to G\times_{S_\rho} S$ yielding an isomorphism of the tangent spaces.
Composing it with (\ref{Tetale1}) we obtain an isomorphism
\begin{equation}\label{Tetale1'}
T_e\, G\times T_\rho\, S\to T_\rho\, {\cal Hom}(\Gamma,G).
\end{equation}
The image of $T_e\ G$ in $T_\rho\, {\cal Hom}(\Gamma,G)$ is $T_\rho\, O_\rho.$ Therefore,
\begin{equation}\label{Tetale1''}
T_\rho S\simeq {\cal Hom}(\Gamma,G)/T_\rho\, O_\rho\simeq H^1(\Gamma,Ad\, \rho).
\end{equation}
Since $\rho$ is a simple point of ${\cal Hom}(\Gamma,G)$ we can assume that $S$ is non-singular at $\rho$ by
\cite[Remark III.1.1]{Lu}. Therefore, there exists an isomorphism
\begin{equation}\label{linearization}
T_{\rho}\, (S//S_\rho)\simeq T_0\, (T_{\rho} S//S_\rho)
\end{equation}
by \cite[Lemma III.1]{Lu} (c.f. \cite[Thm 6.4]{PV}).
This isomorphism combined with (\ref{Tetale1''}) and (\ref{Tetale2}) implies the statement.
\qed

Since $O_\rho$ is always reduced, it is a subvariety of $Hom(\Gamma,G)$ in ${\cal Hom}(\Gamma,G).$ A version of the above argument applied to the $G$ action on $Hom(\Gamma,G)$ yields the following:

\begin{theorem}\label{tangent_X2}
(1) If $\rho$ is completely reducible then the projection $Hom(\Gamma,G)\to X_G(\Gamma)$ induces a linear map
$$T_\rho\, Hom(\Gamma,G)/B^1(\Gamma,Ad\, \rho)\to T_{[\rho]}\, X_G(\Gamma).$$
(2) If $\rho$ is good then this is an isomorphism.\\
(3) If $\rho$ is smooth then
$$T_0\Big(\big(T_\rho\, Hom(\Gamma,G)/B^1(\Gamma,Ad\, \rho)\big)//S_\rho\Big)\simeq T_{[\rho]}\, X_G(\Gamma).$$
\end{theorem}

\begin{corollary}\label{XHom-reduced}
If $\rho$ is good then ${\cal X}_G(\Gamma)$ is reduced at $[\rho]$ if and only if $\rho$ is reduced.
\end{corollary}

\begin{proof}
${\cal X}_G(\Gamma)$ is reduced if and only if $$T_{[\rho]}\, X_G(\Gamma)=T_{[\rho]}\, {\cal X}_G(\Gamma)$$ which holds if and
only if $$T_\rho\, Hom(\Gamma,G)=T_\rho\, {\cal Hom}(\Gamma,G)$$ by Theorems \ref{tangent_X}(2) and \ref{tangent_X2}(2).
\end{proof}

\begin{question}
(1) Is there a natural linear map
$$T_0\, \big(H^1(\Gamma, Ad\, \rho)//S_\rho\big)\to T_{[\rho]} {\cal X}_G(\Gamma)$$
for all completely reducible $\rho$, which is an isomorphism for scheme smooth $\rho$?\\
(2) Is this map an injection or surjection for $\rho$ which are not scheme smooth?
\end{question}

The difficulty in answering this question comes from the following three factors:
\begin{enumerate}
\item \'Etale slices are not unique.
\item The isomorphism (\ref{linearization}) is non-canonical for a given slice $S$.
\item Map $\phi$ does not descend to an isomorphism
$$T_0\, \big(H^1(\Gamma, Ad\, \rho)//S_\rho)\big)\to T_{[\rho]} {\cal X}_G(\Gamma)$$
if $H^1(\Gamma,Ad\, \rho)\to T_0\, \left( H^1(\Gamma,Ad\, \rho)//S_\rho\right)$ is not onto.
\end{enumerate}

The following result illuminates the importance of the requirement of $\rho$ being scheme smooth in Theorem \ref{tangent_X}(2).

\begin{theorem}\label{H^1_iff_reduced}
${\cal X}_G^g(\Gamma)$ is reduced iff
$T_{[\rho]}\, X_G(\Gamma)= H^1(\Gamma,Ad\, \rho)$ for a non-empty Zariski open set of $[\rho]$'s in $X_G^g(\Gamma)$.
\end{theorem}

\begin{proof}
$\Rightarrow$ By Theorem \ref{tangent_X}(2), the equality holds for all $[\rho]\in X_G^g(\Gamma)$.\\
$\Leftarrow$ If ${\cal X}_G^g(\Gamma)$ is not reduced then ${\cal X}_G^g(\Gamma)$ is not reduced at $[\rho]$ for an open set $\Omega$ of $[\rho]$'s in $X_G^g(\Gamma)$.
Since the set of smooth and good representations, $Hom^{sg}(\Gamma,G),$ is the non-singular part of the set of good representations, it is an open dense subset of $Hom^g(\Gamma,G).$
Hence $\Omega \cap Hom^{sg}(\Gamma,G)$ is non-empty and open.

By Corollary \ref{XHom-reduced}, $\rho$'s in $\Omega \cap Hom^{sg}(\Gamma,G)$ are not reduced.
Therefore, $$dim\, T_\rho\, Hom(\Gamma,G)< dim\, T_\rho\, {\cal Hom}(\Gamma,G)=dim\, Z^1(\Gamma,Ad\,\rho)$$
by Theorem \ref{Z^1-thm}. Since
$$dim\, T_{[\rho]}\, X_G(\Gamma) = dim\, T_\rho\, Hom(\Gamma,G)-dim\, B^1(\Gamma,Ad\, \rho)$$
for $\rho\in Hom^{sg}(\Gamma,G),$ by Theorem \ref{tangent_X2}(3),
we get $$dim\, T_{[\rho]}\, X_G(\Gamma) < dim\, H^1(\Gamma,Ad\, \rho)$$ for $\rho\in \Omega \cap Hom^{sg}(\Gamma,G).$
\end{proof}

%
\section{Symplecticity of Character Varieties of Surfaces}
\label{s_char_var_surf}
%

Let $G$ be a reductive group and let $\g$ be its Lie algebra.
A bilinear form \mbox{$B:\g \times \g \to \C$} is $Ad$-invariant if $B(Ad(g)x,Ad(g)y)=B(x,y).$

Let $F$ be a closed orientable surface.
For every representation $\rho:\pi_1(F)\to G$ and every $Ad$-invariant bilinear form
$B:\g \times \g \to \C$, the cup product defines a pairing
\begin{equation}\label{omega}
\omega_B: H^1(F, Ad\, \rho)\times H^1(F, Ad \, \rho) \stackrel{\cup}{\longrightarrow} H^2(F,
Ad\, \rho\otimes Ad\, \rho)\stackrel{B}{\longrightarrow} H^2(F,\C)=\C
\end{equation}
which can be also identified with the pairing
\begin{equation}\label{omega_cap}
H^1(F,Ad\, \rho)\times H_1(F,Ad\, \rho)\stackrel{\cap}{\longrightarrow}
H_0(F,Ad\, \rho\otimes Ad\, \rho)\stackrel{B}{\longrightarrow} H_0(F,\C)=\C
\end{equation}
via the Poincar\'e duality with twisted coefficients, \cite{Sp},
$$\cap [F]: H^n(F,Ad\, \rho)\to H_{2-n}(F,Ad\, \rho),$$
where $[F]\in H_2(F,\C)$ is a fundamental class of $F.$

Let $\Gamma$ be a group and let $(C_*,\p)$ be a chain complex of left $\Z\Gamma$-modules.
Let $M_1,M_2$ be left $\Z\Gamma$-modules. If $B: M_1\times M_2\to \C$ is a $\Z\Gamma$-invariant pairing, i.e. $B(rm_1,rm_2)=B(m_1,m_2)$ for every $r\in \Z\Gamma,$ $m_1\in M_1,$ $m_2\in M_2,$
then, by \cite[V \S 3]{Br}, the cap product induces a pairing
\begin{equation}\label{non-deg}
H^n(Hom_{\Z\Gamma}(C_*,M_2),\delta)\times H_n(M_1\otimes_{\Z\Gamma} C_*,\p)\to M_1\otimes_{\Z\Gamma} M_2
\stackrel{B}{\longrightarrow} \C.
\end{equation}
In the above formula $M_1$ is considered as a right $\Z\Gamma$-module via $m\cdot\gamma=\gamma^{-1}\cdot m.$

\begin{lemma}\label{cap}
If $B$ is a duality pairing, i.e. if $B$ induces an isomorphism $M_1\simeq Hom(M_2,\C),$
then the pairing (\ref{non-deg}) is non-degenerate.
\end{lemma}

\begin{proof}
The cochain complex $(Hom_{\Z\Gamma}(C_*,M_2),\delta)$ can be written as
$$(Hom_{\Z\Gamma}(C_*,Hom(M_1,\C)),\delta)=(Hom_{\Z\Gamma}(C_*\otimes M_1,\C),\delta)=$$
$$(Hom (C_*\otimes_{\Z\Gamma} M_1,\C),\delta)=Hom((C_*\otimes_{\Z\Gamma} M_1,\p),\C).$$
Hence
\begin{equation}\label{chain_iso}
H_n(Hom_{\Z\Gamma}(C_*,M_2),\delta)=H_n\left(Hom((C_*\otimes_{\Z\Gamma} M_1,\p),\C)\right).
\end{equation}
Since $\C$ is a divisible group, $Hom(\ \cdot\ , \C)$ is an exact functor in the category of abelian groups. Hence, (\ref{chain_iso}) becomes
$$H_n(Hom_{\Z\Gamma}(C_*,M_2),\delta)=Hom(H_n(C_*\otimes_{\Z\Gamma} M_1,\p),\C).$$
It is easy to verify that this isomorphism is induced by (\ref{non-deg}).
\end{proof}

If $B$ is symmetric then (\ref{omega}) is skew-symmetric. Therefore, Lemma \ref{cap} implies:

\begin{corollary}\label{cor_symp} (c.f. \cite{Go1})
If $B:\g\times \g\to \C$ is symmetric, $Ad$-invariant, and non-degenerate,
then (\ref{omega}) is a symplectic form on $H^1(F,Ad\, \rho).$
\end{corollary}

If $\g$ is simple then the Killing form is unique among symmetric, $G$-invariant, non-degenerate forms on $\g,$ up to a constant multiple.

Let $F$ be a closed orientable surface of genus $\geq 2$ now.
By Corollary \ref{good-smooth}, the space of conjugacy classes of good representations,
$$X_G^{g}(F)=Hom^{g}(\pi_1(F),G)//G= Hom^{g}(\pi_1(F),G)/G$$
is a complex manifold and, by Proposition \ref{surface-reduced} and Theorem \ref{tangent_X}(2),
 $$T_{[\rho]}\, X_G^{g}(F)=H^1(\pi_1(F),Ad\, \rho).$$

\begin{remark}\label{omega_alg}
$\omega_B$ is an ``algebraic" form on $X^{g}_G(F)$, i.e. it is a global section of the second exterior power of the vector bundle of K\"ahler differentials on $X_G^{g}(F).$
In particular, $\omega_B$ is holomorphic.
\end{remark}

Goldman proves by an argument from gauge theory that for every non-degenerate, symmetric, $Ad$-invariant $B,$ $\omega_B$ is closed, \cite{Go1}. Therefore, $(X_G^{g}(F),\omega_B)$ is a holomorphic symplectic manifold. (Note that our Proposition \ref{surface-reduced} is needed to complete Goldman's construction of $\omega_B$.)

Although $X^g_G(torus)$ is empty for most $G$, $X_G(torus)$ is a singular symplectic manifold as well, \cite{Si-ab}.



\section{$3$-manifolds and Lagrangian Subspaces}
\label{s_lagr}

Let $M$ be an orientable compact $3$-manifold with a connected boundary $F.$
The embedding $\p M\hookrightarrow M$ induces a homomorphism $r:\pi_1(F)\to \pi_1(M)$ and a map \mbox{$r^*: X_G(M)\to X_G(F).$}

It is often believed that the non-singular part of the image of $X_G(M)$ in $X_G^g(F)$ is Lagrangian. We investigate this claim in this section.

\begin{theorem}\label{lag1} Consider the symplectic form $\omega_B$ induced by a symmetric non-degenerate bilinear
Ad-invariant form $B$ on $\g.$\\
(1) Let $\rho:\pi_1(M)\to G$ be such that $\rho r$ is good.
Then the image of $T_{[\rho]}\, {\cal X}_G(M)$ in $T_{[\rho]}\, X_G(F)$ is Lagrangian.\\
(2) The non-singular part of $r^*(X_G(M))$ in $X_G^g(F),$ denoted by us by
$Y_G(M),$ is an isotropic submanifold of $X_G^g(F).$
(In particular, every connected component of $Y_G(M)$ of dimension $\frac{1}{2} dim\, X_G(F)$ is Lagrangian.)\\
(3) If a connected component of $Y_G(M)$ contains an equivalence class of a reduced
$G$-representation of $\pi_1(M)$ then it is a Lagrangian submanifold of $X_G(F).$
\end{theorem}

Although we do not know any example of a $3$-manifold $M$ such that $Y_G(M)$ is not Lagrangian, in light of the above theorem and the results of M. Kapovich mentioned in Section \ref{s_char_var_sch} we believe that such examples exist.

Here is an alternative version of the part (3) above:

\begin{theorem}\label{lag2} Let $X_G^s(M)$ be the set of equivalence classes of scheme smooth $G$-representations of $\pi_1(M).$ Then
$X^g_G(F)\cap r^*(X_G^s(M))$ is an immersed Lagrangian submanifold of $X^g_G(F).$
\end{theorem}

Note that $X_G^s(M)$ may be empty, even if there exist good $G$-representations of $\pi_1(M)$.
Note also that $r^*: X_G^s(M)\to X_G^g(F)$ does not have to be an immersion. Indeed, it is easy to choose an example of a $3$-manifold $M$ satisfying the statement of Proposition \ref{3mfld}(2) whose all $G$-representations
are reduced and, hance, $X_G^s(M)$ has an arbitrarily large dimension.
We prove versions of theorems \ref{lag1} and \ref{lag2} for $F$ torus in \cite{Si-ab}.

Theorems \ref{lag1} and \ref{lag2} have applications to Chern-Simons quantum field theory, \cite{Ba,Fr,Gu,JW,We1,We2,Si3}. It is a $2+1$-topological quantum field theory, associating a Hilbert space $H(F)$ to every closed orientable surface $F$  and a vector $I(M)\in H(F)$ to every $3$-manifold $M$ with $\p M=F.$ Although a mathematically rigorous version of this theory exists for compact groups $G$, \cite{RT,BHMV}, its constructions are combinatorial and algebraic in nature, yielding little information about the relations between $I(M)$ and the topology of $M.$ Furthermore, a rigorous construction of this theory is still missing for complex algebraic groups.

One hopes to achieve a geometric construction of Chern-Simons TQFT for all groups through a procedure of geometric quantization which associates Hilbert spaces $H$ to symplectic manifolds $X$ and vectors in $H$ to Lagrangian subspaces of $X,$ \cite{Sn,Wo}. It is natural to expect that Witten's $I(M)\in H(F)$ is the vector associated with $Y_G(M).$ However, several obstacles exist in this approach, the first one being the question whether $Y_G(M)$ is Lagrangian.



Theorems \ref{lag1} and \ref{lag2} are relevant also to Floer homology theory for $3$-manifolds, c.f. \cite{Cu,BC1,BC2}.
If $F\subset M$ is a surface separating $M$ into $M_1$ and $M_2$ and $Y_G(M_1),Y_G(M_2)$ are Lagrangian submanifolds of $X_G(F)$ then one may consider Floer symplectic homology groups for such a splitting. (A difficulty with this approach lies in the non-compactness of $X_G(F)$.)

For every representation $\rho: \pi_1(M)\to G$ the homomorphism $r: \pi_1(F)\to \pi_1(M)$
induces $r^*: H^1(M,Ad\, \rho)\to H^1(F,Ad\, \rho\, r).$ The proofs of Theorems \ref{lag1} and \ref{lag2} are based on the following:

\begin{theorem}\label{main_H1}
For every $\rho: \pi_1(M)\to G,$ $r^*H^1(M,Ad\, \rho)$ is a Lagrangian subspace of the symplectic space
$(H^1(F,Ad\, \rho\, r),\omega_B)$ with respect to every non-degenerate, $Ad$-invariant, symmetric, bilinear form $B$ on $\g.$
\end{theorem}


In particular, for the trivial representation $\rho: \pi_1(M)\to \C^*$,
Theorem \ref{main_H1} implies the following classical result:

\begin{cor}
For every compact, orientable $3$-manifold with a connected boundary $F$ the image of the map
$r^*: H^1(M,\C)\to H^1(F,\C)$ induced by the embedding $r:F\hookrightarrow M$ is a Lagrangian
subspace of $H^1(F,\C)$ with the symplectic form being the cup product.
\end{cor}

\noindent{\it Proof of Theorem \ref{main_H1}:}
(1) We prove that
$$dim\, r^*H^1(M,Ad\, \rho)=\frac{1}{2}\, dim\, H^1(F,Ad\, \rho\, r)$$
first, by filling in the details of the approach of \cite{Fr}. (This approach was communicated to us by Charlie Frohman.)
By Poincar\'e-Lefschetz duality we have
\begin{equation}\label{diag}
\begin{array}{ccccc}
H_2(M,F,Ad\, \rho) & \stackrel{\p}{\longrightarrow} & H_1(F,Ad\, \rho\, r) & \stackrel{r_*}{\longrightarrow} & H_1(M,Ad\, \rho)\\
\downarrow & & \downarrow \eta & & \downarrow\\
H^1(M,Ad\, \rho) & \stackrel{r^*}{\longrightarrow} & H^1(F, Ad\, \rho\, r) & \stackrel{\delta}{\longrightarrow} &
H^2(M,F,Ad\, \rho),\\
\end{array}
\end{equation}
where all vertical maps are isomorphisms induced by Poincar\'e duality.
By Corollary \ref{cor_symp}, the cap product
$$H^1(F,Ad\, \rho\, r)\times H_1(F,Ad\, \rho\, r)\stackrel{\cap}{\longrightarrow} H_0(F,Ad\, \rho\, r\otimes Ad\, \rho\, r)\stackrel{B}{\longrightarrow} H_0(F,\C)=\C,$$
is non-degenerate.
Similarly,
$$H^1(M,Ad\, \rho)\times H_1(M,Ad\, \rho)\stackrel{\cap}{\longrightarrow} H_0(M,Ad\, \rho\otimes Ad\, \rho)\stackrel{B}{\longrightarrow}\C,$$
$$H^2(M,F,Ad\, \rho)\times H_2(M,F,Ad\, \rho)\stackrel{\cap}{\longrightarrow} H_0(M,F,Ad\, \rho\otimes Ad\, \rho)\stackrel{B}{\longrightarrow} \C$$ are non-degenerate by Lemma \ref{cap}.
Consider the isomorphisms
$$H^1(M,Ad\, \rho)\simeq (H_1(M,Ad\, \rho))^*,\quad H^1(F,Ad\, \rho\, r)\simeq
(H_1(F,Ad\, \rho\, r))^*,$$
$$H^2(M,F,Ad\, \rho)\simeq (H_2(M,F,Ad\, \rho))^*$$ defined by these pairings.
Under these identifications, $r^*$ and $r_*$ and $\p$ and $\delta$ are the duals of each other.
Hence $$rank\, r^*=rank\, r_*=rank\, \delta=dim\, H^1(F,Ad\, \rho\, r)-dim\, Ker\, \delta$$
$$=dim\, H^1(F,Ad\, \rho\, r) - rank\, r^*.$$

(2) It remains to be proven that $r^*H^1(M,Ad\, \rho)$ is an isotropic subspace of $H^1(F,Ad\, \rho).$

The pairing (\ref{omega_cap}) identifies $H_1(F,Ad\, \rho\, r)$ with $H^1(F,Ad\, \rho\, r)^*$.
The isomorphism $\eta^{-1}$ of (\ref{diag}) sends $\alpha\in H^1(F, Ad\, \rho\, r)$ to $\eta^{-1}(\alpha)\in H_1(F, Ad\, \rho\, r)$
which under the above identification is the functional $f_\alpha: H^1(F, Ad\, \rho\, r)\to \C,$ $f_\alpha(\beta)=\omega_B(\alpha,\beta).$

By Lemma \ref{cap}, the pairing
\begin{equation}
H^1(M,Ad\, \rho)\times H_1(M,Ad\, \rho)\stackrel{\cap}{\longrightarrow}
H_0(M,Ad\, \rho\otimes Ad\, \rho)\stackrel{B}{\longrightarrow} H_0(M,\C)=\C
\end{equation}
is non-degenerate.
If we use it to identify $H_1(M,Ad\, \rho)$ with $H^1(M,Ad\, \rho)^*$ then $r_*$ in (\ref{diag}) sends
$f_\alpha$ to $f_\alpha r_*: H^1(M,Ad\, \rho)\to \C$. By commutativity and exactness of (\ref{diag}), $f_\alpha r_*=0$ for every $\alpha\in r^*(H^1(M,Ad\, \rho)).$ In other words, $f_\alpha(\beta)=0$
for every $\alpha,\beta\in r^*(H^1(M,Ad\, \rho)).$
\qed

\begin{remark}\label{tangent_rmk}
Let $G$ be reductive and $\rho: \pi_1(M)\to G$ be such that $\rho\, r: \pi_1(F)\to G$
is good. Then $\rho$ is good as well and, by Theorem \ref{tangent_X}(1) and (2), the following diagram commutes:
$$\begin{array}{ccc}
H^1(M,Ad\, \rho) & \stackrel{r^*}{\longrightarrow} & H^1(F,Ad\, \rho\, r)\\
\downarrow \phi & & \downarrow \phi\\
T_{[\rho]} {\cal X}_G(M) & \stackrel{dr^*}{\longrightarrow} & T_{[\rho r]} {\cal X}_G(F),
\end{array}$$
where $\phi$ is the isomorphism of Theorem \ref{tangent_X}(2).
\end{remark}

{\bf Proof of Theorem \ref{lag1}:}
(1) is a direct consequence of Theorems \ref{tangent_X}, \ref{main_H1} and of the remark above.\\
(2) Since $Y_G(M)$ is smooth, it is enough to show the following statement
\begin{equation}\label{isotrop}
T_{[\rho r]}\, Y_G(M)\subset T_{[\rho r]}\, X_G(F)\quad \text{is isotropic}
\end{equation}
for a dense subset of points $[\rho r]\in Y_G(M)$ (in the complex topology). Let $X_G'(M)\subset X_G(M)$ be the non-singular part of $(r^{*})^{-1}Y_G(M).$
Since $r^*(X_G'(M))$ is dense in $Y_G(F),$ it is enough to show that (\ref{isotrop}) holds for a dense subset of points
in $r^*(X_G'(M)).$ By Sard's theorem the there is a dense set of points in $[\rho r]\in r^*(X_G'(M))$ for which
$$T_{[\rho r]}\, Y_G(M)= dr^* T_{[\rho]}\, X_G(M).$$
Therefore, it is enough to show that $dr^* T_{[\rho]}\, X_G(M)$ is isotropic in $T_{[\rho r]}\, X_G(F).$
This follows from (1) and the fact that $T_{[\rho]}\, X_G(M)$ is a subspace of $T_{[\rho]}\, {\cal X}_G(M).$\\
(3) Let $\rho:\pi_1(M)\to G$ be a reduced, irreducible representation whose conjugacy class
belongs to $Y_G(M).$ By Theorem \ref{tangent_X}(2), both $\phi$'s in the diagram of Remark \ref{tangent_rmk} are isomorphisms. By Theorem \ref{main_H1},
$dim\, T_{[\rho]}\, C= \frac{1}{2} dim\, T_{[\rho]}\, X_G(F).$
Since $C$ and $X_G^g(F)$ are smooth, that equality holds for all points of $C.$
Now the statement follows from (2).
\qed

{\bf Proof of Theorem \ref{lag2}:}
By Proposition \ref{g_open}, $X_G^g(F)\subset X_G(F)$ is open and, hence,
$U=(r^*)^{-1}(X_G^g(F))\cap X_G^s(M)$ is an open subset of $X_G^s(M).$
Since $\phi$'s in the diagram of Remark \ref{tangent_rmk} are isomorphisms, by Theorem \ref{main_H1},
$$dr^*: T_{[\rho]}\, U\to T_{[\rho]}\, X_G^g(F)$$ has constant rank for all $[\rho]$ in $U.$
By the Constant Rank Theorem, for every $[\rho]\in U$ there is a neighborhood $V$ of
$[\rho]$ in $U$ and a neighborhood $W$ of $r^*([\rho])$ in $X_G^g(G)$ such that
$r^*(V)\cap W$ a submanifold of $W.$ Consequently, $r^*(U)=r^*(X_G^s(M))\cap X_G^g(F)$
is an immersed submanifold of $X_G^g(F)$.
\qed

\centerline{Department of Mathematics, 244 Math. Bldg.}
\centerline{University at Buffalo, SUNY}
\centerline{Buffalo, NY 14260, USA}
\centerline{asikora@buffalo.edu}
\end{document}